\documentclass{imsart}

\RequirePackage{amsthm,amsmath,amsfonts,amssymb}
\RequirePackage[authoryear]{natbib}
\RequirePackage[colorlinks,citecolor=blue,urlcolor=blue]{hyperref}
\RequirePackage{graphicx}

\usepackage{latexsym}
\usepackage{float}
\usepackage{xcolor}
\usepackage{url}
\usepackage[normalem]{ulem}
\usepackage{soul}
\usepackage{censor}
\usepackage{enumerate}
\usepackage{epstopdf}
\usepackage{color}
\usepackage{multirow}
\usepackage{array}
\usepackage{multicol}
\usepackage[mathscr]{euscript}

\startlocaldefs

\endlocaldefs

\begin{document}

\begin{frontmatter}
\title{Locally correct confidence intervals for a binomial proportion: A new criteria for an interval estimator}
\runtitle{Locally correct confidence intervals for a binomial proportion}

\begin{aug}
\author{\fnms{Paul H.} \snm{Garthwaite}},
\author{\fnms{Maha W.} \snm{Moustafa}}
\and
\author{\fnms{Fadlalla G.} \snm{Elfadaly}}
\address{School of Mathematics and Statistics, The Open University, Milton Keynes, MK7 6AA, UK.}
\runauthor{Garthwaite et al.}
\end{aug}

\begin{abstract}
Well-recommended methods of forming `confidence intervals' for a binomial proportion give interval estimates that do not actually meet the definition of a confidence interval, in that their coverages are sometimes lower than the nominal confidence level. The methods are favoured because their intervals have a shorter average length than the Clopper-Pearson (gold-standard) method, whose intervals really are confidence intervals. Comparison of such methods is tricky -- the best method should perhaps be the one that gives the shortest intervals (on average), but when is the coverage of a method so poor that it should not be classed as a means of forming confidence intervals?

As the definition of a confidence interval is not being adhered to, another criterion for forming interval estimates for a binomial proportion is needed. In this paper we suggest a new criterion; methods which meet the criterion are said to yield {\em locally correct confidence intervals}. We propose a method that yields such intervals and prove that its intervals have a shorter average length than those of any other method that meets the criterion. Compared with the Clopper-Pearson method, the proposed method gives intervals with an appreciably smaller average length. The mid-$p$ method also satisfies the new criterion and has its own optimality property.
\end{abstract}



\end{frontmatter}


\section{Introduction}
A good number of methods of forming a confidence interval for a binomial success parameter ($p$) have been studied. For example, Volsett (1993) compares thirteen methods, Newcombe (1998) compares seven methods and Brown et al. (2001) examine eleven. These studies only consider methods that aim to form equal-tail confidence intervals, which can be constructed from one-sided intervals. That is, if $(l, \, u)$ is a $(1- 2\alpha)$ confidence interval for $p$, then $(l, \, 1)$ and $(0, \, u)$ are one-sided $(1-\alpha)$ confidence intervals for $p$. The `gold-standard' method for forming equal-tail confidence intervals is the Clopper-Pearson method (Leemis and Trivedi, 1996; Jovanovic and Levy, 1997). This interval estimator meets the definition of a method for forming equal-tail confidence intervals: the coverage of its one-sided $(1-\alpha)$ confidence intervals is {\em guaranteed} to be at least $(1-\alpha)$ for any value of $p$. Moreover, among those methods that meet the definition, its intervals have the shortest average length.

The drawback of the Clopper-Pearson method is that its coverage is above the nominal confidence level for almost all values of $p$, and for modest sample sizes this conservatism is noticeable. This is illustrated in the left-hand diagram in Figure~1, which plots the coverage of its 97.5\% upper one-sided intervals against $p$ when samples are drawn from a bin$(20,p)$ distribution. The saw-tooth pattern in the plot results from the discreteness of the sample space and arises with any method of forming confidence intervals for a binomial proportion. It can be seen that coverage is commonly above 98.5\% and sometimes exceeds 99.5\%. One approach to reducing conservatism is to construct confidence intervals that are not constrained to have equal tails. Crow (1956), Blyth and Still (1983) and Casella (1986) adopt this approach and produce effective methods that reduce conservatism while maintaining the nominal confidence level. In practice, however, equal-tail intervals seem to be much preferred to two-sided intervals with unequal tails. Here we only consider one-sided intervals, and two-sided intervals with equal tails.

\begin{figure}
 \includegraphics[clip=true, viewport=0 30 1290 485,scale=0.38]{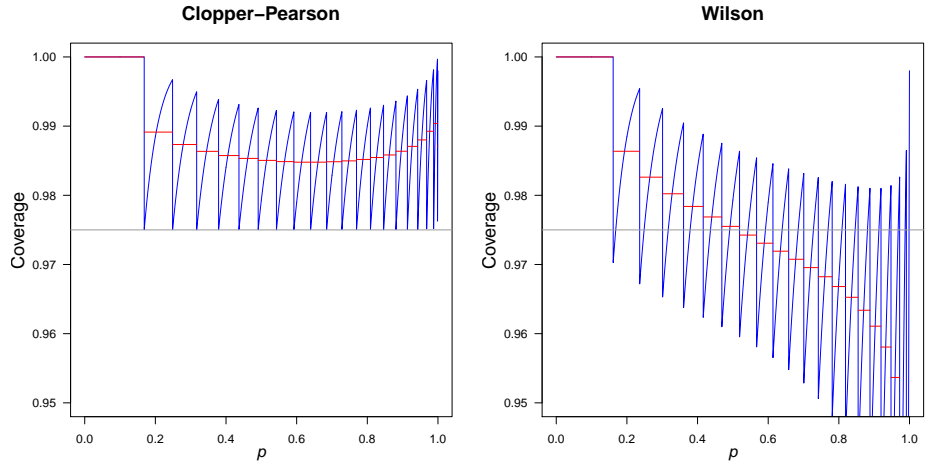}
\caption{Coverage of upper one-sided 97.5\% confidence intervals for the Clopper-Pearson and Wilson methods for a fixed sample size ($n$) of 20 and variable $p$. Short horizontal (red) lines show the average coverage between consecutive spikes}
\end{figure}

Agresti and Coull (1998, page 119) argue that the conservatism of the Clopper-Pearson method makes it ``\ldots inappropriate to treat this approach as optimal for statistical practice''. Instead, they advocate interval estimators that give shorter intervals for which the coverage probability is usually quite close to the nominal level, though coverage is less than that level for some values of $p$. An interval estimator that is often recommended under this criteria is the Wilson 
method (Wilson, 1927), whose coverage is plotted in the right-hand diagram in Figure~1 (again for 97.5\% upper one-sided intervals for sampling from a bin$(20,p)$ distribution). For most values of $p$ its coverage is quite close to the nominal level and, because it is sometimes liberal and generally less conservative than the Clopper-Pearson method, on average its intervals are shorter.

This pragmatic attitude underlies almost all work on methods of forming equal-tail confidence intervals for a binomial proportions -- interval estimators  are examined that do not strictly give confidence intervals (their coverage is sometimes too small) but their coverages are reasonably close to the nominal level. Recommendations as to which interval estimator is ``best'' are made on the basis of coverage and the average length of their intervals. Simplicity of the methods is sometimes considered as well. However, methods giving shorter intervals tend to have lower coverage and an appropriate trade-off between accuracy of coverage and average interval length is partly a matter of opinion
and recommendations as to which method is best cannot be clear-cut.

This situation is regrettable. Assuming equal-tails intervals are required, it seems to the writers that an interval estimator for a binomial proportion should not actually aim to produce confidence intervals.  If that were the aim, then the Clopper-Pearson method would be the method of choice because it gives equal-tails intervals with a shorter average length than any other method that genuinely yields confidence intervals. Nor should the notion of a confidence interval be relaxed in a non-rigorous fashion -- simply requiring coverage to be `close' to the nominal level creates ambiguity as to whether one interval estimator is better than another.  Instead, we feel that interval estimators should be required to meet some new criterion. With an appropriate criterion, there should be some interval estimators that

\begin{enumerate}
\item[(a)] satisfy the new criterion,
\item[(b)] give intuitively sensible intervals,
\item[(c)] give intervals whose average length is acceptably short.
\end{enumerate}

Under such a criterion, it would be reasonable to restrict attention to interval estimators with the above properties and the best estimator would be the one with the shortest average interval length. The challenge is to find an appropriate new criterion.

In this paper we propose a criterion that reflects the saw-tooth pattern of coverage that is illustrated in Figure~1. For an upper one-sided interval, the coverage increases monotonically as $p$ increases from one spike to the next, until it drops when the next spike is reached. Similarly, for a lower one-sided interval, the coverage decreases monotonically between spikes as $p$ increases. Our new criterion is that, for one-sided intervals, the average coverage between any pair of consecutive spikes must never be less than the nominal level. Interval estimators that meet this criterion will be said to yield {\em locally correct confidence} (LCC) {\em  intervals}. In Figure~1, short horizontal (red) lines show the average coverage between pairs of consecutive spikes. The Clopper-Pearson intervals are LCC intervals while those given by Wilson's method are not.

Two-sided equal-tail $(1-2\alpha)$ LCC intervals are formed by constructing lower-tail and upper-tail $(1-\alpha)$ LCC intervals.  We continue to refer to the nominal level as the confidence level. Most importantly when, for example, we specify an interval and say it is a 95\% LCC interval, the statement would be accurate rather than just approximately true.

In Section~2 we define the criterion more precisely and in Section~3 we present a novel interval estimator that yields LCC intervals. We refer to the new estimator as the {\em optimal locally correct} (OLC) {\em  method} and prove that it yields intervals with a smaller average length than any other interval estimator that yields LCC intervals. We also examine whether intervals given by the OLC method seem sensible and identify properties of the new estimator that have been proposed in the literature as being desirable.

The proofs of most results require separate consideration of each binomial sample size that is of interest. The conservatism of the Clopper-Pearson method dissipates as the sample size increases, becoming negligible, so there is little reason to seek an alternative to the gold- standard Clopper-Pearson method for large sample sizes. Hence, asymptotic results are of limited relevance and instead we prove results for all sample sizes up to 200.

In Section~4 we compare the OLC method with a number of methods that have been recommended for forming equal-tail confidence intervals. One of these methods is the mid-$p$ method. We find that this method yields LCC intervals if $\alpha$ is less than 0.1, which includes all the values of $\alpha$ that are commonly used in forming confidence intervals. We also identify a property of the mid-$p$ method that seems to have been overlooked in the past. The mid-$p$ method is well-recommended as a means of forming confidence intervals (Mehta and Walsh, 1992; Vollset, 1993; Berry and Armitage, 1995; Agresti and Gottard, 2005), suggesting that its intervals are acceptably short. As the average lengths of its intervals are a little larger than with the new method, this implies that the new method also gives intervals that are acceptably short. To facilitate use of the OLC method, in an appendix we provide tables of its two-sided 95\% and 99\% intervals for sample sizes up to 30. Intervals for other sample sizes and confidence levels can be obtained through a Shiny R application available at \url{https://olcbinomialci.shinyapps.io/binomial/}.

Requiring interval estimators to yield LCC intervals is a criterion that meets requirements (a)--(c) listed above. There are, no doubt, alternative criteria that also meet these requirements. However, some seemingly plausible criteria fail to meet them. In Section~5, two such criteria are described that we briefly considered but they give intervals that could not be considered sensible. Concluding comments are given in Section~6.

\section{Locally correct confidence intervals}
We first consider upper-tail (upper one-sided) intervals. Let $X$ denote a binomial variate based on $n$ trials with success probability $p$ and suppose an interval estimator gives $(0,u_x)$ as its upper-tail estimate for $p$ when $x$ is the observed value of $X$. For a sensible estimator,
\begin{equation} \label{eq1}
0 < u_0  < u_1  < \ldots  < u_n  \leq 1,
\end{equation}
and we assume that (\ref{eq1}) holds. The coverage probability of the interval estimator depends upon the value of $p$ and is the probability that the random interval $(0,u_X)$ contains $p$. We denote this coverage probability by $C_u(p)$. When $u_{i-1} < p \leq u_i$,
\begin{equation} \label{eq2}
C_u(p) \, = \,P(X \geq i \, | \, p)\, = \, \sum_{x=i}^n \binom{n}{x}p^x (1-p)^{n-x}.
\end{equation}

For $x=0,1, \ldots,n$, the difference between $C_u(u_x)$ and $C_u(u_x+\delta)$ does not tend to 0 as $\delta \rightarrow 0$, but equals $\binom{n}{x}u_x^x (1-u_x)^{n-x}$. This is the reason that the coverages in Figure~1 have a saw-tooth appearance. The points of the teeth (the spikes) occur where $p$ equals $u_0, \ldots, u_n$ and the coverage drops by $\binom{n}{x}u_x^x (1-u_x)^{n-x}$ at $p=u_x$.

For a good interval estimator, how should $C_u(p)$ vary with $p$\,? If the nominal confidence level is $(1-\alpha)$, then $C_u(p)$ should exceed $(1-\alpha)$ when $p$ is just before a spike, as $C_u(p)$ follows a cycle with its largest values just before spikes. If the estimator is not to be very conservative, then $C_u(p)$ should be less than $(1-\alpha)$ when $p$ is just after a spike, as then $C_u(p)$ is at the lowest part of its cycle. However, the extent to which the estimator is liberal should be restricted. It seems reasonable to require the {\em average coverage within each cycle} to exceed or equal the target confidence level. That is the restriction we impose. A cycle -- the interval between two spikes -- is quite small. Hence, while the coverage need not equal (or exceed) the confidence level at individual values of $p$, it must do so on average over quite narrow ranges of $p$. We say that an interval estimator that meets this requirement gives {\em locally correct confidence} (LCC) {\em intervals}.\vspace{.08in}\\
{\em Definition 1}.\, Suppose an interval estimator gives $(0,u_x)$ as its
upper-tail interval for $p$ when $X=x$, that $u_0, \ldots, u_n$ satisfy (1) and that $u_n=1$. If, for $i=1,\ldots,n$,
\begin{equation} \label{eq3}
\frac{1}{u_i-u_{i-1}} \,\int_{p=u_{i -1}}^{u_i} C_u(p) \, dp \, \geq 1- \alpha
\end{equation}
then the interval estimator gives {\em upper-tail LCC intervals} with confidence level $1-\alpha$.\vspace{.08in}\\
In Definition 1 it is assumed that $u_n=1$; otherwise the average coverage over the interval $(u_n,1)$ would be 0, which is inconsistent with the required coverage in other intervals. Definition 2 is the corresponding definition for lower-tail intervals.\vspace{.08in}\\
{\em Definition 2}.\, Suppose an interval estimator gives $(l_x,1)$ as its lower-tail interval for $p$ when $X=x$ and that $0=l_0 < l_1 < \ldots < l_n<1$. Define the coverage probability, $C_l(p)$, by
\begin{equation} \label{eqB}
C_l(p) \, = \,P(X \leq i \, | \, p)\, = \, \sum_{x=0}^i \binom{n}{x}p^x (1-p)^{n-x}
\end{equation}
for $l_{i} < p \leq l_{i+1}$. If, for $i=0,\ldots,n-1$,
\begin{equation} \label{eq5}
\frac{1}{l_{i+1}-l_{i}} \,\int_{p=l_{i}}^{l_{i+1}}
C_l(p) \, dp \, \geq 1-\alpha
\end{equation}
then the interval estimator gives {\em lower-tail LCC intervals} with confidence level $1-\alpha$.\vspace{.08in}

Two-sided equal-tail LCC interval estimators are defined in terms of one-sided LCC intervals.\vspace{.08in}\\
{\em Definition 3}.\, Suppose that, for $(x=0,\ldots,n)$, an interval estimator gives $(l_x,u_x)$ as its two-sided equal-tail intervals for $p$ when $X=x$. Then it gives {\em equal-tail LCC intervals} with confidence level $(1-2\alpha)$ if and only if $ \{ (l_x,1), x=0,\ldots,n \} $ and  $ \{ (0,u_x), x=0,\ldots,n \} $ are sets of one-sided lower-tail and upper-tail LCC intervals, respectively, each with confidence level $(1-\alpha)$.\vspace{.08in}\\
An interval estimator that gives equal-tail LCC intervals will be referred to as an {\em LCC interval estimator}.

\section{A new interval estimator}
We propose an interval estimator that uses a straightforward iterative algorithm to obtain one-sided interval estimates. The following proposition underpins the algorithm. Proofs of propositions are given in supplementary material.\vspace{.1in}\\
{\bf Proposition 1.}\, Suppose $1 \leq n \leq 200$ and $0.0001 <\alpha <0.27 $. Suppose also that
\begin{equation} \label{eq7}
\frac{1}{u_i-u_{i-1}} \,\int_{p=u_{i -1}}^{u_i}\sum_{x=i}^n \binom{n}{x}p^x (1-p)^{n-x} \, dp \, = 1- \alpha
\end{equation}
and $u_i >u_{i-1}$ for $i=j+1, j+2, \ldots, n$; $j=0, \ldots, n-1$. Then there is a unique $u_{j-1}$ such that $u_j >u_{j-1}>0$ and equation~(\ref{eq7}) holds when $i=j$. Also, each $u_j$ is a monotonic decreasing function of $\alpha$ for $0.0001 < \alpha <0.27$.\vspace{.15in}\\
A necessary condition in Proposition~1 is that $\alpha<0.27$. For larger values of $\alpha$, sometimes $P(X \geq i |p=u_i)$ is less than $ 1- \alpha$ and then equation (\ref{eq7}) cannot be satisfied with $u_{i-1} \leq u_i$. However, $\alpha$ is generally set equal to $0.05$, $0.025$ or $0.005$ for a confidence interval, and for an inter-quartile range $\alpha=0.25$. Hence Proposition~1 covers the values of $n$ and $\alpha$ of practical interest. The monotonicity property given in Proposition~1 is used in proofs of further propositions.

For an upper-tail interval with confidence level $1-\alpha$ and sample size $n$, the algorithm sets $u_n$ equal to 1 and then sequentially determines $u_{n-1}, u_{n-2}, \ldots , u_0$. Given $u_i$, a simple numerical search is used to find $u_{i-1}$ that satisfies equation~(\ref{eq7}). From Proposition~1, there is always a unique $u_{i-1}$ for which equation~(\ref{eq7}) holds. The new interval estimator sets $(0,u_i)$ as the upper-tail interval when $X=i$. From Proposition~1, we also have that $0 <u_0 < \cdots < u_n =1$. Consequently, under definition 1 the new estimator gives upper-tail LCC intervals.

The new interval estimator uses similar steps to form lower-tail intervals. The iterative search starts by putting $l_0 =0$ and then $l_{1}, \ldots , l_n$ are determined sequentially. Given $l_i$, the value $l_{i+1}$ is found that satisfies
\begin{equation} \label{eq8}
\frac{1}{l_{i+1}-l_{i}} \,\int_{p=l_{i}}^{l_{i+1}}
\sum_{x=0}^i \binom{n}{x}p^x (1-p)^{n-x} \, dp \, = 1-\alpha
\end{equation}
for $i=0, \ldots, n-1$. Then $(l_i,1)$ is the $1-\alpha$ lower-tail interval when $X=i$  and, under definition 2, the estimator gives lower-tail LCC intervals.
Two-sided intervals are obtained by combining the endpoints of one-sided intervals. Thus, when $X=i$ the new estimator gives $(l_i, u_i)$ as the two-sided equal-tails interval for a confidence level of $1-2\alpha$, where $(l_i, 1)$ and $(0, u_i)$ are the $1-\alpha$ lower-tail  and upper-tail LCC intervals, respectively. From its construction, the new estimator is an LCC interval estimator.

Examples showing the coverage of the estimator are given in Figure~2, which plots coverage against $p$ for upper-tail intervals with nominal confidence levels of 97.5\% and 99.5\% for sample sizes 8, 20 and 50. The positions of spikes have similar characteristics to their positions in Figure~1: apart from the gap from 0 to the first spike, the spikes seem fairly evenly spaced and any trend in the size of the inter-spike interval is smooth. However, the coverages in Figure~2 are evenly spread around the nominal level, unlike the conservative coverage of the Clopper-Pearson method (left-hand diagram in Figure~1). Also, the coverages in Figure~2 have no long-term trend, unlike the coverages with Wilson's method (right-hand diagram in Figure~1), which are mostly above the nominal level for small values of $p$ and below the nominal level for large values of $p$. Consequently, Figure~2 suggests that the new method gives sensible interval estimates. By design, with the new method the average coverage between two spikes will always equal the nominal level.

An important feature of an interval estimator is the length of its intervals. If an estimator gives $(l_x, u_x)$ as its interval estimate when $X=x$, then the expected length of its interval, $L_{n}(p)$ say, is given by
\begin{equation} \label{e11}
L_{n}(p)= \sum_{x=0}^n (u_x-l_x)\binom{n}{x}p^x (1-p)^{n-x}
\end{equation}
and, averaging over $p$, its average expected length (AEL) is
\begin{equation} \label{eq12}
E [ L_{n}(p) ] = \int_{p=0}^{1} L_{n}(p) \, dp .
\end{equation}
This definition holds for two-tail intervals that do not necessarily have equal tails and also for one-tail intervals -- in equation~(\ref{e11}) $l_x$ is set equal to 0 for upper-tail intervals and $u_x$ is set equal to 1 for lower-tail intervals.

We will refer to our new method of forming interval estimates as the {\em optimal locally correct} (OLC) method because it has the optimality properties given in Proposition~2, whose proof is given in the supplementary material.\vspace{.08in}\\ \newpage
{\bf Proposition 2.} \, Suppose $1 \leq n \leq 200$ and $0.0001< \alpha < 0.27 $. Then
\begin{enumerate}
\item[(i)] for a one-tail interval with confidence level $1-\alpha$, the OLC method has the smallest AEL of any interval estimator that gives one-tail LCC intervals;
\item[(ii)] for an equal-tail interval with confidence level $1-2\alpha$, the OLC method has the smallest AEL of any LCC estimator.
\end{enumerate}

\begin{figure}
\includegraphics[clip=true, viewport=0 10 1250 500,scale=0.7]{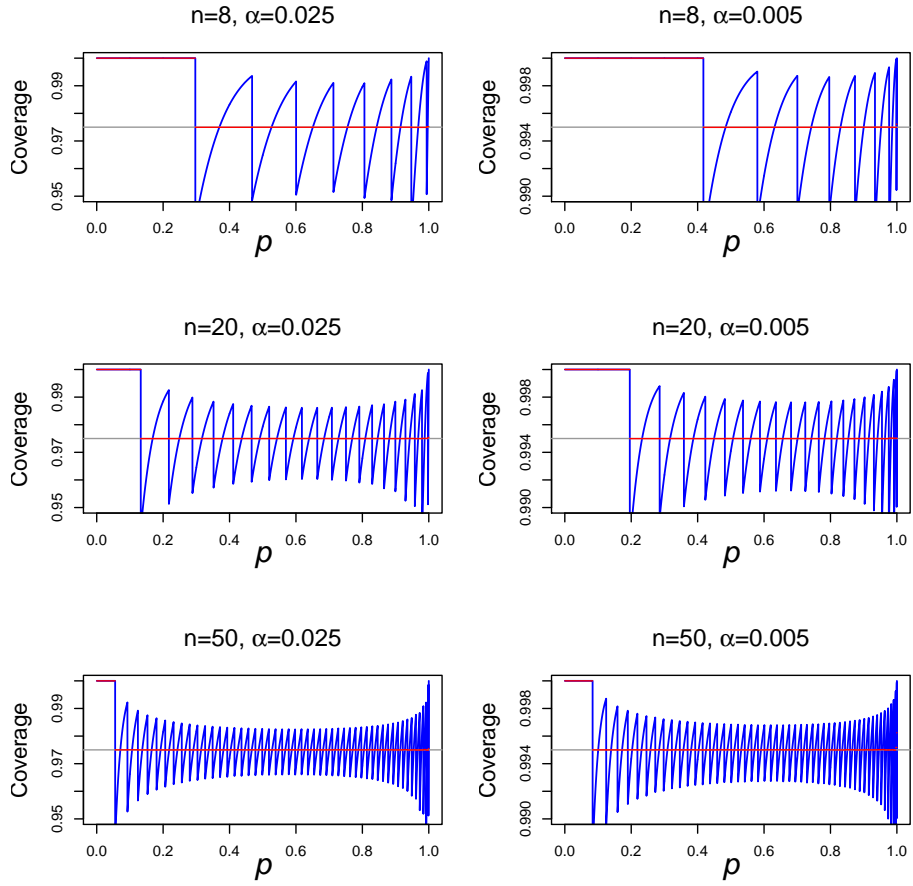}
\caption{Coverage of upper-tail LCC intervals given by the new estimator for sample sizes of 8, 20 and 50, and nominal confidence levels of 97.5\% and 99.5\%}
\end{figure}

The following are desirable properties in a method of forming confidence regions [see, for example, Blyth and Still (1983) and Schilling and Doi (2014)].
\begin{enumerate}
\item [1.]{\em Interval valued.}\, A confidence region should be an interval and not a collection of disjoint intervals.
\end{enumerate}
The remaining properties assume that the confidence region is a two-tail interval. When $X=x$, the sample size is $n$ and the confidence level is $1-2\alpha$, denote this interval as $(L(x,n, \alpha), \, U(x,n, \alpha))$.
\begin{enumerate}
\item [2.] {\em Equivariance.} \, As the binomial distribution is invariant under the transformation $X \rightarrow n-X$; $p \rightarrow 1-p$, confidence intervals should also be invariant under these transformations. That is, $ L(x,n, \alpha)$ should equal $1-U(n-x,n, \alpha)$ for $x=0, \ldots, n$.
\item [3.] {\em Nesting.} \, If two confidence intervals have different confidence levels then, for any given $x$ and $n$, the interval for the higher confidence level should contain the interval for the lower confidence level. If the confidence levels are $1-\alpha_1$ and $1-\alpha_2$ with $\alpha_1 < \alpha_2$, this requires $(L(x,n, \alpha_2), \, U(x,n, \alpha_2)) \in (L(x,n, \alpha_1), \, U(x,n, \alpha_1))$.
\item [4.] {\em Monotonicity in $x$.} \, For fixed $n$ and $\alpha$, the endpoints should be increasing in $x$. This requires $L(x+1,n, \alpha) > L(x,n, \alpha)$ and $U(x+1,n, \alpha) > U(x,n, \alpha)$.
\item [5.] {\em Monotonicity in $n$.} \, For fixed $x$ and $\alpha$, the lower endpoint should be non-increasing in $n$ and the upper endpoint should be decreasing in $n$. This requires $L(x,n+1, \alpha) \leq L(x,n, \alpha)$ and $U(x,n+1, \alpha) < U(x,n, \alpha)$.
\end{enumerate}
For the property of monotonicity in $n$, the ``greater or equal'' inequality cannot be a ``strictly greater'' inequality for the lower limit because, regardless of $n$, the lower limit should be 0 when $x=0$. Instead, the property implies that, when an additional trial results in a failure, the limits of the confidence interval should not increase. In conjunction with the invariance property it also implies that, when an additional trial results in a success, the limits of the confidence interval should not decrease.

Suppose $1 \leq n \leq 200$ and $0.0001 <\alpha < 0.27 $, so that Proposition~1 applies. Then clearly the OLC method will have the first four properties, as (i) its estimate is an interval and never a set of disjoint intervals, (ii) equivalent procedures are used to construct lower and upper tails, (iii) the nesting property follows from the monotonicity property given by Proposition~1 and (iv) interval endpoints increase as $x$ increases. Proposition~3, given below, shows that the OLC method also gives interval endpoints that have the `monotonicity in $n$' property. Hence it seems clear that the OLC method meets the above requirements for being a well-behaved interval estimator.\vspace{.08in}\\
{\bf Proposition 3.} \, The endpoints of intervals given by the OLC method satisfy $L(x,n+1, \alpha) \leq L(x,n, \alpha)$ and $U(x,n+1, \alpha) < U(x,n, \alpha)$ for $1 \leq n \leq 200$ and $0.0001< \alpha < 0.27 $.

\section{Comparison with other methods}
In this section we compare the OLC method with the following six methods of forming interval estimates: Clopper-Pearson, mid-$p$, Wald, Agresti-Coull,  Wilson, and Jeffreys methods. An optimality property of the mid-$p$ method is also given.

\subsection{Clopper-Pearson and  mid-$p$ methods}
To simplify notation, we let $(l_i, \, u_i)$ denote a method's $1-2\alpha$ confidence interval when $i$ successes occur in a sample of size $n$. The Clopper-Pearson method determines $(l_i, \, u_i)$ by inverting equal-tail tests of the hypothesis $H_0 : \, p=p_0$. Thus $l_i$ satisfies $P(X \geq i \, | \, p=l_i)=  \alpha$ and $u_i$ satisfies $P(X \leq i \, | \, p=u_i)= \alpha$,
except that $l_0 =0$ and $u_n =1$.
Of the methods we consider, this is the only one that strictly meets the definition of a method for forming confidence intervals. Since it gives equal-tail confidence intervals, it also gives LCC intervals. As noted earlier, its conservative coverage (cf. Figure~1) leads to intervals that are frequently considered unnecessarily long.

The mid-$p$ method reduces the conservatism of the Clopper-Pearson method by halving  the probability of the observed result (Newcombe, 1998). That is,  $l_i$ satisfies
\begin{equation} \label{eq13}
P(X > i \, | \, p=l_i) +
 {\textstyle\frac{1}{2}} P(X = i \, | \, p=l_i) = \alpha
\end{equation}
and $u_i$ satisfies
\begin{equation} \label{eq14}
P(X < i \, | \, p=u_i)
+  {\textstyle\frac{1}{2}} P(X = i \, | \, p=u_i) = \alpha
\end{equation}
except, by definition, $l_0 =0$ and $u_n =1$, as otherwise the coverage is 0 for $p<l_0$ and $p>u_n$ (Agresti and Gottard, 2005). While the mid-$p$ method does not give interval estimates that meet the definition of confidence intervals, the method does give LCC intervals for $n<200$  and $0.0001<\alpha<0.1$. The range of $\alpha$ ($0.0001<\alpha<0.1$) is more restrictive than with the OLC method and, in particular, does not include interquartile ranges, but it does include all levels of confidence that are commonly of interest when forming confidence intervals. The result is given in the following proposition.\vspace{.08in}\\
{\bf Proposition 4.} \, For $1 \leq n \leq 200$ and $\alpha \in (0.0001, \,0.1 ) $, the mid-$p$ method gives LCC intervals.\vspace{.1in}\\
The coverage of 97.5\% upper-tail intervals for the mid-$p$ method for a bin$(20, p)$ distribution is plotted against $p$ in the upper-left panel of Figure~3. The  spikes in the plot are spaced fairly regularly and the actual coverage always crosses the nominal coverage level between consecutive spikes.

\begin{figure}
\includegraphics[clip=true, viewport=0 20 1250 500,scale=0.7]{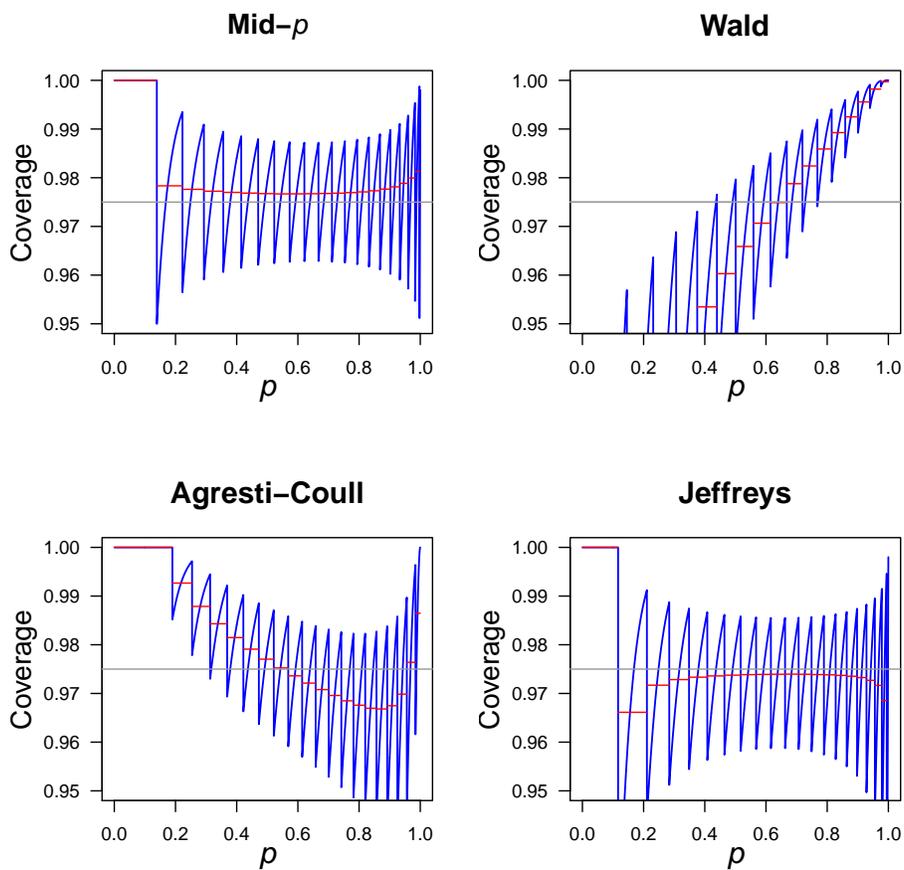}
\caption{Coverage of upper one-sided 97.5\% confidence intervals for the mid-$p$, Wald, Agresti-Coull and Jeffreys methods for a fixed sample size of 20 and variable $p$}
\end{figure}

Agresti and Gottard (2007) compare the coverages of the mid-$p$ and Clopper-Pearson methods as $p$ varies uniformly between 0 and 1. They conclude that ``the mid-$p$ approach is an excellent one to adopt if one hopes to achieve close to the nominal level in using a method repeatedly for various studies in which $p$ itself varies''
(Agresti and Gottard, 2007,  p.6455).
In fact, among an important class of methods of forming confidence intervals, the mid-$p$ method gives one-tail confidence intervals whose coverage is optimally close to the nominal confidence level for any value of $p$. This property is given in Proposition~5. It seems unlikely that the property has not been noted before, but we have been unable to find a reference to it in the literature, so a proof of the proposition is given in the supplementary material. We only give the property for upper-tail intervals but the corresponding result obviously holds for lower-tail intervals.\vspace{.08in}\\
{\bf Proposition 5.} \,
Consider the class of methods of forming upper-tail confidence intervals that (i) do not involve randomisation (i.e. confidence interval are determined by $x$ and $n$, and do not involve the value of a further hypothetical random variable), and (ii) satisfy $0 \leq u_0 \leq \ldots \leq u_n=1$. Then the absolute error in the coverage probability, $|C_u(p) - (1-\alpha)|$, is as small or smaller for the mid-$p$ method as for any method in the class, for any value of $p$.

The class of methods in Proposition~5 includes all sensible methods of forming equal-tail confidence intervals that do not involve randomisation so, in particular, it includes all methods discussed in this section.  The property is quite strong because it relates to every value of $p$, and hence gives other properties. For example, in comparing methods of interval estimation, it is common to examine the average absolute error in coverage or the root mean-square error in coverage, where averaging is over $p \sim U(0,1)$. Under either of these measures, Proposition~5 implies that the mid-$p$ method would be the optimal method of forming one-tail confidence intervals.

\subsection{Wald, Agresti-Coull,  Wilson, and Jeffreys methods}
The remaining methods we consider do not give LCC intervals.
The Wald method is the simple method that is commonly taught in introductory statistics courses. Its interval endpoints are given by
\begin{equation}  \label{eq15}
l_i =\widehat{p} -z \{\widehat{p}(1-\widehat{p})/n \} ^{1/2}\, \, \mbox{ and } \, \,
u_i =\widehat{p} +z \{\widehat{p}(1-\widehat{p})/n \} ^{1/2} ,
\end{equation}
where $\widehat{p} = i/n$ and $z$ is is the $1-\alpha$ quantile of the standard normal distribution. Wald confidence intervals have been criticized heavily in, for example, Brown et al. (2001) and Vos and Hudson (2008). One fault is that their intervals tend to have low coverage and Agresti and Coull (1998) suggested a simple adjustment that appreciably improves the coverage of the 95\% confidence intervals. The adjustment, which gives the Agresti-Coull method, is to add two ``successes'' and two ``failures'' to the sample, so the endpoints of its confidence interval are $\tilde{p} \, \pm \, z \{\tilde{p}(1-\tilde{p})/(n+4) \} ^{1/2} $, where $\tilde{p} = (i+2)/(n+4)$. The coverages of Wald and Agresti-Coull 97.5\% upper-tail intervals for the bin$(20, p)$ distribution are plotted against $p$ in the upper-right and lower-left panels of Figure~3. In contrast to the mid-$p$ coverage, the coverages show clear trends.  The coverage of Wald intervals tends to be liberal for small values of $p$ and conservative for large values, while the trend is in the opposite direction for Agresti-Coull intervals, and much less pronounced. Agresti-Coull intervals are recommended by Brown et al. (2001) as a simple method of forming confidence intervals for a binomial proportion when the sample size exceeds 40.

When the length of confidence intervals is of paramount concern, perhaps the most well-recommended methods are the Wilson method and Jeffreys method (see, for example, Brown et al., 2001). Wilson intervals are based on inversion of the score test, so they are also known as score intervals. The interval endpoints are
\[\frac {n\hat{p}+z^{2}/2}{n+z^{2}} \pm \frac{z}{n+z^{2}} \left\{ n\hat{p}(1-\hat{p})+ \frac{z^{2}}{4} \right\} ^{1/2}. \]
A plot of the coverage of 97.5\% upper-tail intervals against $p$ for bin$(20,p)$ distributions was given in the right panel of Figure~1. Coverage is a little conservative for small values of $p$ and quite liberal for values of $p$ near 1. These characteristics are reversed for lower-tail intervals, and so the conservatism and liberalism partly balance out when the coverage of two-tail intervals is plotted against $p$. (It is more common to plot the coverage of two-tail intervals, rather than one-tail intervals, when the coverage of Wilson's method looks better.)

Jeffreys method is a Bayesian approach that takes a Beta$(\frac{1}{2}, \, \frac{1}{2})$ distribution as the prior distribution, which is Jeffreys' choice of noninformative prior distribution for sampling from a binomial model. The sample consists of $x$ successes in $n$ trials and leads to Beta$(x+\frac{1}{2}, \, n-x+\frac{1}{2})$ as the posterior distribution for $p$. Jeffreys interval estimate is the $1-2\alpha$ equal-tail credible interval given by this posterior distribution,  except for setting $l_0=0$ and $u_n =1$ (Brown et al., 2001). Thus $l_i$ and $u_i$ are set equal to the $\alpha$ and $1-\alpha$ quantiles of the Beta$(i+\frac{1}{2}, \, n-i+\frac{1}{2})$ distribution for $i=1,\ldots, n-1$. Brown et al. (2001) note that Jeffreys intervals are always within Clopper-Pearson intervals and are approximately equal to mid-$p$ intervals. The coverage of 97.5\% upper-tail intervals for Jeffreys method is given in the lower-right panel of Figure~3. Its coverage resembles that of the mid-$p$ method and does not display the trend found with Wilson's method.

\subsection{Coverage and length of intervals}
We will restrict attention to upper-tail intervals and two-tail intervals. As defined in Section~2, $C_u(p)$ is the coverage of an upper-tail interval estimator and is the probability that the random interval $(0,u_X)$ contains $p$. We define the quantity $T_u$ as
\begin{equation}  \label{eq16}
T_u = \frac{1}{1-u_0}\int_{p=u_0}^1 C_u(p) \, dp,
\end{equation}
and refer to it as the {\em truncated average coverage}. In calculating $T_u$, values of $p$ in the range $(0,u_0)$ are excluded because within that range the coverage is 1, so the coverage for $p \leq u_0$ differs radically from the coverage for $p>u_0$. Consequently, it is more informative to give the values of both $T_u$ and $u_0$. Other average coverages may be calculated from these: the average coverage of a one-tail interval over the full range $(0,1)$ equals $ \{ (1-u_0 )T_u +u_0 \} $ and the average coverage for two-tail intervals over the range $(0,1)$ equals $ 2 \{ (1-u_0 )T_u +u_0 \} -1$.

Table~1 gives both the values $T_u$ and $u_0$ for the methods described earlier, for each of $\alpha = 0.05$, 0.025 and 0.005, and $n = 8$, 20 and 50. The results indicate the merits of the OLC method. From its construction, the OLC method necessarily has a truncated average cover that equals the nominal $\alpha$ (apart from rounding error). In contrast, the Clopper-Pearson method is very conservative, mid-$p$ is slightly conservative, and Agresti-Coull is generally conservative, while the Wilson, Wald and Jeffreys methods are consistently liberal. (The Wald method is particularly liberal with values of $T_u$ that are almost always far below the nominal values.)

\begin{table}

\caption{ \label{table:tab1} Truncated average coverage ($T_u$) of upper-tail $1-\alpha$ intervals and smallest upper limit $(u_0)$ of seven methods of forming interval estimates, for $\alpha$ = 0.05, 0.025, 0.005 and sample sizes $(n)$ of 8, 20 and 50.}
\centering
\begin{tabular}{ cccccccccc } \hline
	$\alpha $ & $n$ & statistic &	Clopper	&	Mid-$p$	&	Agresti	&	 Wilson	&	 Wald	 &	 Jeff.	&	OLC	\vspace{-.04in}\\
	 &  &  & Pearson	&	&	Coull	&	&	 &	&	\\ \hline
0.05	&	8	&	$T_u$	&	0.976	&	0.956	&	0.949	&	 0.941	&	 0.852	 &	 0.941	&	0.950	\vspace{-.05in}\\
0.05	&	8	&$u_0$	&	0.312	&	0.250	&	0.293	&	0.253	 &	0.000	 &	 0.208	 &	0.239	\vspace{.05in}\\
0.05	&	20	&	$T_u$	&	0.971	&	0.954	&	0.953	&	 0.947	&	 0.903	 &	 0.946	&	0.950	\vspace{-.05in}\\
0.05	&	20	&$u_0$	&	0.139	&	0.109	&	0.141	&	0.119	 &	0.000	 &	 0.091	 &	0.105	\vspace{.05in}\\
0.05	&	50	&	$T_u$	&	0.966	&	0.952	&	0.953	&	 0.949	&	 0.928	 &	 0.948	&	0.950	\vspace{-.05in}\\
0.05	&	50	&$u_0$	&	0.058	&	0.045	&	0.062	&	0.051	 &	0.000	 &	 0.038	 &	0.043	\vspace{.05in}\\
0.025	&	8	&	$T_u$	&	0.989	&	0.979	&	0.972	&	 0.966	&	 0.867	 &	 0.969	&	0.975	\vspace{-.05in}\\
0.025	&	8	&$u_0$	&	0.369	&	0.312	&	0.372	&	0.324	 &	0.000	 &	 0.262	 &	0.297	\vspace{.05in}\\
0.025	&	20	&	$T_u$	&	0.986	&	0.977	&	0.976	&	 0.972	&	 0.923	 &	 0.972	&	0.975	\vspace{-.05in}\\
0.025	&	20	&$u_0$	&	0.168	&	0.139	&	0.190	&	0.161	 &	0.000	 &	 0.117	 &	0.133	\vspace{.05in}\\
0.025	&	50	&	$T_u$	&	0.983	&	0.976	&	0.977	&	 0.974	&	 0.950	 &	 0.974	&	0.975	\vspace{-.05in}\\
0.025	&	50	&$u_0$	&	0.071	&	0.058	&	0.085	&	0.071	 &	0.000	 &	 0.049	 &	0.056	\vspace{.05in}\\
0.005	&	8	&	$T_u$	&	0.998	&	0.996	&	0.991	&	 0.988	&	 0.882	 &	 0.993	&	0.995	\vspace{-.05in}\\
0.005	&	8	&$u_0$	&	0.484	&	0.438	&	0.509	&	0.453	 &	0.000	 &	 0.379	 &	0.417	\vspace{.05in}\\
0.005	&	20	&	$T_u$	&	0.998	&	0.996	&	0.994	&	 0.992	&	 0.941	 &	 0.994	&	0.995	\vspace{-.05in}\\
0.005	&	20	&	$u_0$	&	0.233	&	0.206	&	0.289	&	 0.249	&	 0.000	 &	 0.177	&	0.196	\vspace{.05in}\\
0.005	&	50	&	$T_u$	&	0.997	&	0.995	&	0.995	&	 0.994	&	 0.970	 &	 0.995	&	0.995	\vspace{-.05in}\\
0.005	&	50	&	$u_0$	&	0.101	&	0.088	&	0.139	&	 0.117	&	 0.000	 &	 0.075 & 0.084 \vspace{.01in}\\ \hline
\end{tabular}
\end{table}

A small value of $u_0$ is desirable, as then the range over which the coverage equals 1 is small. On that basis the Wald method does exceptionally well, as $u_0$ always equals 0 for that method. However, the coverage of the Wald method is too liberal for it to be the preferred method of forming interval estimates. Based on $u_0$, the OLC method is a little poorer than Jeffreys (a consistently liberal method), but a little better than mid-$p$, and much better than Clopper-Pearson, Agresti-Coull  and Wilson.

As noted in the introduction, it has been argued that a good interval estimator should (i) give short intervals, and (ii) give coverage probabilities that are usually quite close to the nominal level (see, for example, Agresti and Coull (1998)). To examine the latter criterion, the root mean-square error (RMSE) of each method's coverage was determined over the truncated range $(u_0,1)$. This RMSE is given by
\begin{equation} \label{eq17a}
\mbox{RMSE} = \left[\frac{1}{1-u_0} \int_{u_0} ^1 \left\{ C_u(p) - (1-\alpha) \right\} ^2 \, dp \right] ^{1/2},
\end{equation}
where $1-\alpha $ is the nominal confidence level. The RMSE for each method is given in Table~2 for $\alpha=0.05,$ 0.025, 0.005 and $n=8,$ 20, 50. From Proposition~5, the mid-$p$ method has the minimum possible RMSE of any method of forming interval estimates that does not use randomisation. Consequently, the mid-$p$ method has the smallest RMSE in every row of Table~2. The new method, OLC, has the second smallest RMSE in every row and is always only a little poorer than the mid-$p$ (its RMSE is never more than 20\% bigger). The RMSE of the Clopper-Pearson was sometimes more than 45\% bigger than the RMSE of the mid-$p$ method and each of the other methods has an RMSE that is at least 80\% bigger for some combination of $\alpha$ and $n$, with the Wald method often doing extremely badly. Hence OLC, while not the optimal method, has a very respectable RMSE.

\begin{table}

\caption{ \label{table:tab2} Root mean-square error (RMSE) of coverage of upper-tail $1-\alpha$ intervals for seven methods of forming interval estimates, for $\alpha$ = 0.05, 0.025, 0.005 and sample sizes $(n)$ of 8, 20 and 50.}
\centering

\begin{tabular}{ ccccccccc } \hline
	$\alpha $ & $n$ & 	Clopper	&	Mid-$p$	&	Agresti	&	Wilson	&	 Wald	 &	 Jeff.	 &	 OLC	\vspace{-.04in}\\
	 &  & 	Pearson	&		&	Coull	&	& &		 &	\\ \hline
			0.05	&	8	&	0.0290	&	0.0224	&	0.0267	&	 0.0313	&	 0.2249	&	 0.0314	&	0.0242	\vspace{-.05in}\\
			0.05	&	20	&	0.0235	&	0.0166	&	0.0188	&	 0.0213	&	 0.1481	&	 0.0212	&	0.0175	\vspace{-.05in}\\
			0.05	&	50	&	0.0180	&	0.0119	&	0.0139	&	 0.0151	&	 0.0972	&	 0.0144	&	0.0124	\vspace{.05in}\\
			0.025	&	8	&	0.0153	&	0.0121	&	0.0177	&	 0.0237	&	 0.2334	&	 0.0186	&	0.0134	\vspace{-.05in}\\
			0.025	&	20	&	0.0125	&  0.0091	&   0.0118  &	 0.0157 &    0.1522	&    0.0123	&   0.0097  \vspace{-.05in}\\
			0.025	&	50	&	0.0097	&	0.0066	&	0.0088	&	 0.0110	&	 0.0983	&	 0.0083	&	0.0069	\vspace{.05in}\\		
			0.005	&	8	&	0.0033	&	0.0027	&	0.0085	&	 0.0147	&	 0.2407	&	 0.0051	&	0.0032	\vspace{-.05in}\\		
			0.005	&	20	&	0.0028	&	0.0021	&	0.0041	&	 0.0088	&	 0.1557	&	 0.0032	&	0.0023	\vspace{-.05in}\\		
			0.005	&	50	&	0.0022	&	0.0015	&	0.0028	&	 0.0057	&	 0.0997	&	 0.0021	&	0.0017	\vspace{.01in}\\	\hline
			
\end{tabular}
\end{table}

Turning to the length of intervals, two-tail intervals are examined as the length of one-tail intervals varies too much with $p$\,: the length of one-tail intervals is approximately proportional to $p$ for upper-tail intervals and to $1-p$ for lower-tail intervals. In Figure~4, the expected lengths of 95\% two-tail intervals are plotted against $p$ for sample sizes 8, 20 and 50. To aid comparison, the expected length for the (new) OLC method is included in all plots. For all values of $p$ and each combination of $n$ and $\alpha$, the expected lengths of the OLC, mid-$p$, Agresti-Coull, Wilson and Jeffreys intervals are all very similar, and a little smaller than the expected lengths of the Clopper-Pearson intervals. Wald intervals have a much smaller expected length than other methods when $p$ is quite large or quite small, but it only achieves this by giving coverages that are well-short of the nominal confidence level (c.f. Table~1).

\begin{figure}
\includegraphics[clip=true, viewport=-10 10 1250 650,scale=0.7]{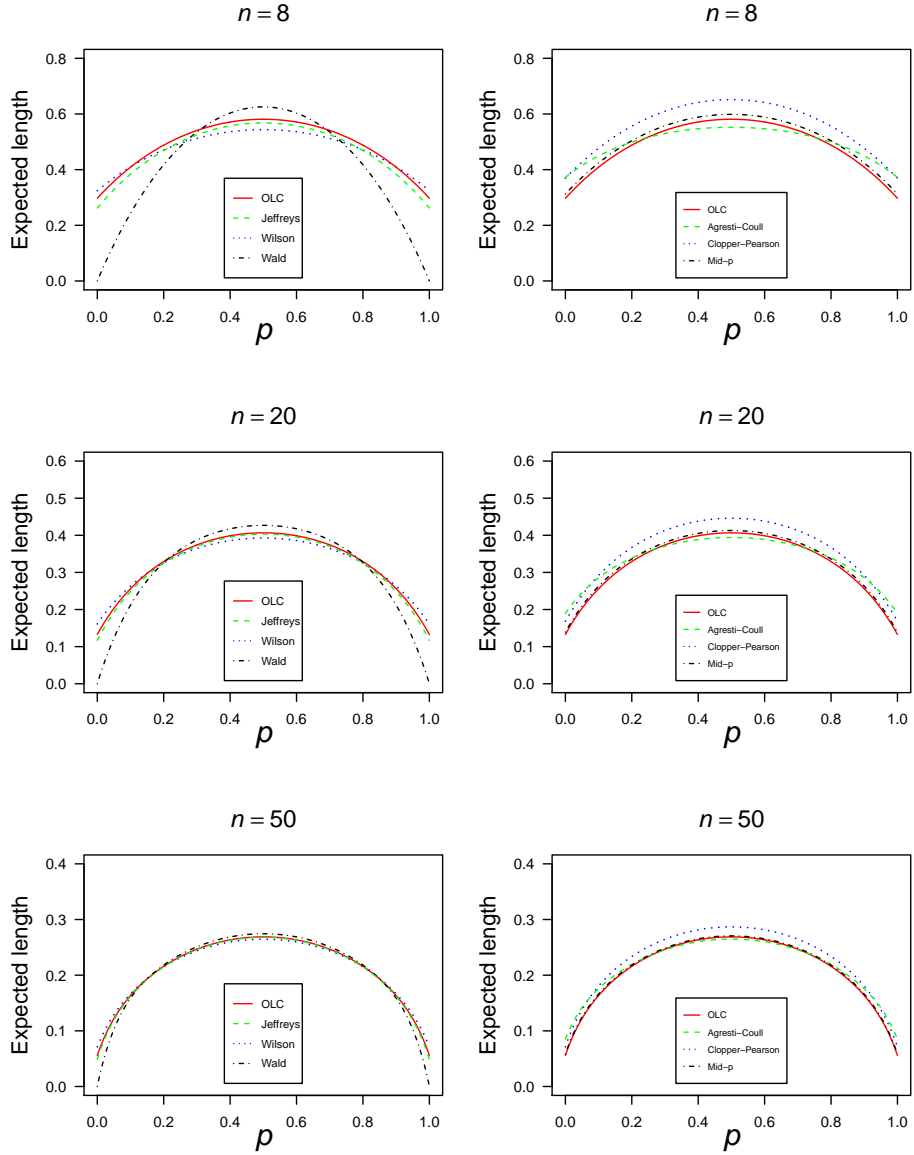}
\noindent
\caption{Expected lengths of two-sided 95\% interval estimates for the OLC, Jeffreys, Wilson and Wald methods (left-hand panels) and  the OLC, Agresti-Coull, Clopper-Pearson and mid-$p$ methods (right-hand panels), plotted against $p$ for sample sizes of 8, 20 and 50}
\end{figure}

Equation~(\ref{eq12}) defines the average expected length (AEL) of intervals given by a method. Important questions about the AEL of the new method are:\vspace{-.1in}
\begin{itemize}
\item [1.]The OLC method gives intervals that have a longer AEL than some of the other methods considered here. But those methods do not give intervals with locally correct coverage. The question is whether the difference in AEL is small, so that there is little cost in requiring intervals to have locally correct coverage.
\item[2.] From Proposition~2, the OLC method gives intervals with a shorter AEL than the AEL of the Clopper-Pearson method. However, the gold-standard Clopper-Pearson method gives intervals that meet the definition of a confidence interval. The question is whether the OLC method reduces the AEL by enough to compensate for failing to meet that definition.
\end{itemize}
To answer those questions requires a means of deciding when a difference is ``sufficiently large'' and when it is ``sufficiently small''. Now, in the literature it has been variously argued that the Wilson, Agresti-Coull, mid-$p$ and Jeffreys methods should be favoured over the Clopper-Pearson method because their intervals are shorter. Hence, the differences in AEL between these methods and the Clopper-Pearson method are considered sufficiently large to compensate for failing to meet the definition of a confidence interval. This provides a measure of a ``sufficiently large'' difference.

\begin{table}
\caption{\label{table:tab3} Average expected length (AEL) of two-tail $1-2\alpha$ intervals for seven methods of forming interval estimates, for  $\alpha$ = 0.05, 0.025, 0.005 and $n$ = 8, 20, 50.}
\centering
\begin{tabular}{ ccccccccc } \hline
	$\alpha $ & \multicolumn{1}{c}{$n$} & 	Clopper	&	Mid-$p$	&	 Agresti	&	 Wilson	 &	 Wald	 &	 Jeff.	 &	OLC	\vspace{-.04in}\\
	 &  & 	Pearson	&		&	Coull	&	& &		 &	\\ \hline
0.050 & 8 & 0.497 & 0.435 & 0.427 & 0.407 & 0.372 & 0.402 & 0.421\vspace{-.05in}\\
0.050 & 20 & 0.317 & 0.283 & 0.284 & 0.275 & 0.268 & 0.273 & 0.278\vspace{-.05in}\\
0.050 & 50 & 0.197 & 0.181 & 0.182 & 0.179 & 0.178 & 0.178 & 0.179\vspace{.05in}\\
0.025 & 8 & 0.561 & 0.508 & 0.499 & 0.474 & 0.427 & 0.472 & 0.492\vspace{-.05in}\\
0.025 & 20 & 0.366 & 0.335 & 0.337 & 0.325 & 0.316 & 0.323 & 0.328\vspace{-.05in}\\
0.025 & 50 & 0.231 & 0.215 & 0.218 & 0.213 & 0.211 & 0.212 & 0.213\vspace{.05in}\\
0.005 & 8 & 0.673 & 0.634 & 0.614 & 0.586 & 0.520 & 0.597 & 0.617\vspace{-.05in}\\
0.005 & 20 & 0.457 & 0.431 & 0.435 & 0.417 & 0.403 & 0.417 & 0.423\vspace{-.05in}\\
0.005 & 50 & 0.295 & 0.281 & 0.286 & 0.278 & 0.275 & 0.276 & 0.278\vspace{-.01in}\\ \hline
	\end{tabular}
\end{table}

Table~3 gives the AEL of the methods considered earlier for each combination of $\alpha = 0.05$, 0.025, 0.005 and $n = 8$, 20, 50. The Wald intervals have the smallest or equally smallest AEL for each combination but, as noted earlier (c.f. Table~1), this is obtained by having an actual coverage that is well below the nominal coverage. Apart from the Wald method and Clopper-Pearson, the AEL of the new method is usually similar in size to that of the other methods, and is always shorter than the mid-$p$ method. Hence, in response to question 1, it is reasonable to conclude that the AEL of the new method compares satisfactorily with that of other methods, so there is little cost in requiring intervals to be locally correct. Also, the new method improves upon the AEL of Clopper-Pearson by amounts of similar size to the improvements given by the other methods. As these other methods are commonly favoured over the Clopper-Pearson method, the response to the second question is also positive -- it is fair to say that the new method improves on the AEL of Clopper-Pearson by enough to compensate for failing to meet the definition of a confidence interval.

\section{Other criteria for defining an interval estimator}
Requiring intervals to be locally correct is just one criterion that might be used in an alternative definition of a ``confidence interval''. In this section we briefly examine two other criteria that seem plausible alternatives but find that they lead to optimal estimators that are often unsatisfactory. We continue to assume that equal-tail intervals are required.

A simple criterion is to require the average coverage of a one-tail interval to be no less than the nominal level. In the notation of equations~(\ref{eq2}) and (\ref{eqB}), this criterion requires $\int_0 ^1 C_u(p) \, dp \geq 1-\alpha$ for upper-tail intervals and $\int_0 ^1 C_l(p) \, dp \geq 1-\alpha$ for lower-tail intervals. The optimal estimator would minimize the average expected length of two-tail intervals under this criterion. Unfortunately, examples show that it will often give endpoints that are identical for a number of different $x$ values. To illustrate, when $n=10$ the optimal estimator gives $(0.3445, \, 0.6555)$ as the 90\% equal-tail interval for $p$ if $x$ equals 3, and it gives exactly the same interval if $x$ equals, 4, 5, 6 or 7. This is clearly unsatisfactory. (Intervals given by the optimal estimator were calculated using the Rsolnp package in R (Ghalanos and Theussl, 2015).)

A more stringent criterion is to partition the $[0, \,1]$ interval into a number of subintervals of equal length and require the average coverage of a one-tail interval to equal or exceed the nominal coverage in each sub-interval. Again, the optimal estimator would minimize the average expected length of two-tail intervals. This criteria offers some flexibility in its application, as the number of subintervals given by the partition must be chosen. However, if the number of subintervals is small, the same problem can arise as with previous criterion -- different values of $x$ will sometimes lead to identical upper endpoints or identical lower endpoints. Also, if the number is large, then the coverage suffers from the form of conservatism that besets the Clopper-Pearson method. Indeed, as the number of subintervals increases, the optimal intervals under this criterion become indistinguishable from those given by the Clopper-Pearson method. Moreover, there is not always an in-between ground that suffers from neither of these problems. As an example, when $n=20$ and $\alpha=0.025$, at least one pair of $x$ values give the same upper-endpoint to the interval estimate unless the number of subintervals is at least 48. Figure~5 shows the coverage of the optimal estimator for an upper-tail 97.5\% interval when the number of subintervals is 48. The short horizontal lines give the average coverage as $p$ varies across each subinterval. The average coverage is equal to the nominal level in fewer than half the subintervals and is often substantially above it, indicating  marked conservatism.

\section{Concluding comments}
This paper aimed to find a satisfactory criterion for choosing an interval estimator for a binomial proportion. As noted in Section~2, with an appropriate criterion there should be some interval estimators that (a) satisfy the new criterion; (b) give intuitively sensible intervals; and (c) give intervals whose average length is acceptably short.

\begin{figure}
\includegraphics[clip=true, viewport=0 0 500 500,scale=0.6]{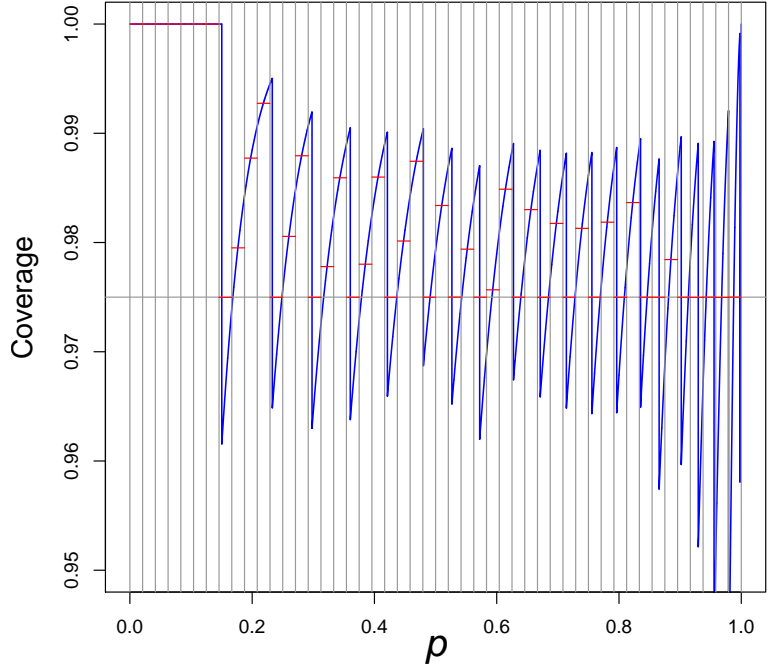}
\caption{Coverage of upper one-sided 97.5\% confidence intervals for $n=20$ and 47 subintervals when the average coverage in each subinterval must be no less than the nominal value. Short horizontal (red) lines show the average coverage in each subinterval}
\end{figure}

As a criterion we proposed that an interval estimator should yield locally correct confidence intervals, meaning that the average coverage between any pair of consecutive spikes should be no smaller than the nominal confidence level. Three of the methods that were examined met this criterion: the Clopper-Pearson method, the mid-$p$ method and the new OLC method. The Clopper-Pearson method arguably does not satisfy point (c), as over the years the conservative length of its intervals has motivated the construction of many other methods of forming confidence intervals for a binomial proportion. One of these other methods is the mid-$p$ method, which has been advocated because of its shorter intervals (Vollset, 1993; Agresti and Gottard, 1995). The mid-$p$ method gives intuitively sensible intervals so it meets point (a)--(c) above.

Turning to the OLC method, in the examples given in Figure~2 the method gave end-points that are fairly evenly spaced with coverages that are balanced around the nominal confidence level. This has also been the case in every other example we have examined and, in particular, there has never been a hint of the defects found with the methods considered in Section~5. Blyth and Still (1983) and Schilling and Doi (2014) list some properties that are desirable in an interval estimator, such as equivariance and monotonicity. We showed that the OLC method has these for the many combinations of $n$ and $\alpha $ that were examined through extensive computation [$1 \leq n \leq 200$ and $\alpha$ in $ ( 0.0001, 0.27 ) $]. Hence it is reasonable to conclude that the new method gives sensible intervals and meets point (b). Regarding the third point, length of intervals, six methods of constructing equal-tail confidence intervals were compared in Section~4.3. The six methods include the mid-$p$, Wilson, Agresti-Coull and Jeffreys methods, which have each been recommended in preference to the Clopper-Pearson method because of the lengths of their intervals. The intervals given by the OLC method had an average expected length that was shorter than Clopper-Pearson intervals and comparable to those given by other methods for all values of $n$ and $\alpha$ that were examined, except when a method gave intervals whose coverage was decidedly smaller than the nominal confidence level. Hence the OLC appears to give intervals that are acceptably short, and meets point (c).

As there are at least two methods that meet points (a)--(c), it can be concluded that requiring intervals to be locally correct is a reasonable criterion to place on an interval estimator. (The results in Section~5 show that finding a suitable criterion is a non-trivial task.) Choosing between the mid-$p$ method and the OLC method is tricky because each has an optimality property. On the one hand, for any value of $p$, the coverage of one-tail intervals is as close to the nominal level as possible when intervals are determined using the mid-$p$ method (c.f. Proposition~5). On the other hand, the average expected length of intervals is smaller with the OLC method than with any other method that gives locally correct intervals (c.f. Proposition~2). However, in choosing between estimators, both coverage and length of intervals are important. Hence either (i)  a restriction should be placed on coverage and methods should be differentiated on the basis of interval width, or (ii) a restriction should be placed on interval width and methods differentiated on the basis of coverage. Requiring an interval to be locally correct places a restriction on coverage, so differentiating between methods on the basis of interval length seems appropriate, in which case the OLC method is the preferred interval estimator.

\subsection*{Package}
A Shiny R application is available at \url{https://olcbinomialci.shinyapps.io/binomial/} that determines one-tail or two-tail intervals for any values of $n$, $x$ and $\alpha<0.5$.


\subsection*{Acknowledgements}
%
%
The second author was supported by a studentship from The Open University, UK.


\subsection*{Supplementary material}
Proofs of Propositions~1--5 are given in supplementary material.

\newpage
\begin{appendix}
\section*{Appendix: Tables of lower and upper endpoints of 95\% and 99\% equal-tail confidence intervals given by the OLC method}

\noindent
Table A1:  Limits of 95\% equal-tail OLC interval when $x$ successes are observed in a sample of size $n$

\scriptsize

\noindent
\begin{tabular}{cccccccccccccc} \hline
$x$	&&	\multicolumn{2}{c}{ $n=1$}		&	\multicolumn{2}{c}{$n=2$}			 &	 \multicolumn{2}{c}{$n=3$}			&	 \multicolumn{2}{c}{$n=4$}			 &	 \multicolumn{2}{c}{$n=5$}			&	\multicolumn{2}{c}{$n=6$}			 \\  \hline					 
0	&&	0.0000	&	0.9500	&	0.0000	&	0.7396	&	0.0000	&	 0.6038	&	 0.0000	 &	 0.5003	&	0.0000	&	0.4291	&	0.0000	&	 0.3735	 \\					
1	&&	0.0500	&	1.0000	&	0.0252	&	0.9748	&	0.0169	&	 0.8444	&	 0.0127	 &	 0.7400	&	0.0101	&	0.6470	&	0.0085	&	 0.5759	 \\					
2	&&	    	&	    	&	    	&	    	&	0.1556	&	 0.9831	&	 0.1115	 &	 0.8885	&	0.0870	&	0.8049	&	0.0714	&	 0.7239	 \\					
3	&&	    	&	    	&	     	&	    	&	     	&	    	 &	    	 &	    	 &	0.1951	&	0.9130	&	0.1566	&	 0.8434	 \\ &&&&&&&&&&&&&\\ \hline					
$x$	&&	\multicolumn{2}{c}{$n=7$}			&	 \multicolumn{2}{c}{$n=8$}			 &	 \multicolumn{2}{c}{$n=9$}			 &	 \multicolumn{2}{c}{$n=10$}			&	 \multicolumn{2}{c}{$n=11$}			 &	\multicolumn{2}{c}{$n=12$}			\\ \hline					 
0	&&	0.0000	&	0.3314	&	0.0000	&	0.2971	&	0.0000	&	 0.2696	&	 0.0000	 &	 0.2465	&	0.0000	&	0.2271	&	0.0000	&	 0.2105	 \\					
1	&&	0.0072	&	0.5160	&	0.0063	&	0.4680	&	0.0056	&	 0.4272	&	 0.0051	 &	 0.3933	&	0.0046	&	0.3639	&	0.0042	&	 0.3388	 \\					
2	&&	0.0606	&	0.6584	&	0.0526	&	0.6006	&	0.0464	&	 0.5525	&	 0.0416	 &	 0.5105	&	0.0377	&	0.4746	&	0.0344	&	 0.4430	 \\					
3	&&	0.1309	&	0.7726	&	0.1125	&	0.7129	&	0.0986	&	 0.6585	&	 0.0878	 &	 0.6121	&	0.0792	&	0.5706	&	0.0721	&	 0.5345	 \\					
4	&&	0.2274	&	0.8691	&	0.1937	&	0.8063	&	0.1688	&	 0.7519	&	 0.1496	 &	 0.7011	&	0.1344	&	0.6569	&	0.1220	&	 0.6167	 \\					
5	&&		   	&	    	&	    	&	    	&	0.2481	&	 0.8312	&	 0.2186	 &	 0.7814	&	0.1954	&	0.7340	&	0.1768	&	 0.6921	 \\					
6	&&	    	&	    	&	    	&	    	&	    	&	    	 &	    	 &	    	 &	0.2660	&	0.8046	&	0.2398	&	 0.7602	 \\ &&&&&&&&&&&&&\\ \hline					
$x$	&&	\multicolumn{2}{c}{$n=13$}			&	 \multicolumn{2}{c}{$n=14$}			 &	 \multicolumn{2}{c}{$n=15$}			 &	 \multicolumn{2}{c}{$n=16$}			&	 \multicolumn{2}{c}{$n=17$}			 &	\multicolumn{2}{c}{$n=18$}			\\ \hline					 
0	&&	0.0000	&	0.1962	&	0.0000	&	0.1836	&	0.0000	&	 0.1726	&	 0.0000	 &	 0.1628	&	0.0000	&	0.1541	&	0.0000	&	 0.1462	 \\					
1	&&	0.0039	&	0.3167	&	0.0036	&	0.2974	&	0.0034	&	 0.2802	&	 0.0032	 &	 0.2650	&	0.0030	&	0.2512	&	0.0028	&	 0.2389	 \\					
2	&&	0.0317	&	0.4154	&	0.0294	&	0.3908	&	0.0273	&	 0.3691	&	 0.0256	 &	 0.3495	&	0.0241	&	0.3319	&	0.0227	&	 0.3159	 \\					
3	&&	0.0662	&	0.5022	&	0.0611	&	0.4737	&	0.0568	&	 0.4479	&	 0.0531	 &	 0.4249	&	0.0498	&	0.4040	&	0.0469	&	 0.3851	 \\					
4	&&	0.1117	&	0.5812	&	0.1030	&	0.5490	&	0.0956	&	 0.5203	&	 0.0892	 &	 0.4941	&	0.0836	&	0.4705	&	0.0786	&	 0.4489	 \\					
5	&&	0.1614	&	0.6534	&	0.1486	&	0.6189	&	0.1376	&	 0.5873	&	 0.1281	 &	 0.5588	&	0.1199	&	0.5326	&	0.1127	&	 0.5087	 \\					
6	&&	0.2184	&	0.7206	&	0.2005	&	0.6835	&	0.1854	&	 0.6501	&	 0.1724	 &	 0.6192	&	0.1611	&	0.5912	&	0.1512	&	 0.5652	 \\					
7	&&	0.2794	&	0.7816	&	0.2559	&	0.7441	&	0.2361	&	 0.7086	&	 0.2192	 &	 0.6765	&	0.2046	&	0.6464	&	0.1918	&	 0.6189	 \\					
8	&&	    	&	    	&	    	&	    	&	0.2914	&	 0.7639	&	 0.2700	 &	 0.7300	&	0.2516	&	0.6990	&	0.2356	&	 0.6698	 \\					
9	&&	    	&	    	&	    	&	    	&	    	&	    	 &	    	 &	    	 &	0.3010	&	0.7484	&	0.2815	&	 0.7185	 \\ &&&&&&&&&&&&&\\ \hline					
$x$	&&	\multicolumn{2}{c}{$n=19$}			&	 \multicolumn{2}{c}{$n=20$}			 &	 \multicolumn{2}{c}{$n=21$}			 &	 \multicolumn{2}{c}{$n=22$}			&	 \multicolumn{2}{c}{$n=23$}			 &	\multicolumn{2}{c}{$n=24$}			\\ \hline					 
0	&&	0.0000	&	0.1391	&	0.0000	&	0.1327	&	0.0000	&	 0.1268	&	 0.0000	 &	 0.1215	&	0.0000	&	0.1165	&	0.0000	&	 0.1120	 \\					
1	&&	0.0027	&	0.2277	&	0.0025	&	0.2174	&	0.0024	&	 0.2081	&	 0.0023	 &	 0.1995	&	0.0022	&	0.1916	&	0.0021	&	 0.1843	 \\					
2	&&	0.0215	&	0.3015	&	0.0204	&	0.2882	&	0.0194	&	 0.2761	&	 0.0185	 &	 0.2649	&	0.0177	&	0.2547	&	0.0169	&	 0.2451	 \\					
3	&&	0.0443	&	0.3678	&	0.0420	&	0.3520	&	0.0399	&	 0.3374	&	 0.0380	 &	 0.3240	&	0.0363	&	0.3116	&	0.0348	&	 0.3002	 \\					
4	&&	0.0742	&	0.4292	&	0.0703	&	0.4110	&	0.0668	&	 0.3944	&	 0.0636	 &	 0.3790	&	0.0607	&	0.3647	&	0.0580	&	 0.3515	 \\					
5	&&	0.1063	&	0.4868	&	0.1006	&	0.4667	&	0.0954	&	 0.4480	&	 0.0908	 &	 0.4309	&	0.0866	&	0.4149	&	0.0828	&	 0.4000	 \\					
6	&&	0.1425	&	0.5414	&	0.1348	&	0.5194	&	0.1278	&	 0.4991	&	 0.1215	 &	 0.4802	&	0.1158	&	0.4627	&	0.1107	&	 0.4464	 \\					
7	&&	0.1805	&	0.5934	&	0.1706	&	0.5698	&	0.1616	&	 0.5479	&	 0.1536	 &	 0.5276	&	0.1463	&	0.5086	&	0.1397	&	 0.4910	 \\					
8	&&	0.2216	&	0.6430	&	0.2091	&	0.6179	&	0.1980	&	 0.5947	&	 0.1880	 &	 0.5730	&	0.1790	&	0.5528	&	0.1708	&	 0.5338	 \\					
9	&&	0.2644	&	0.6903	&	0.2493	&	0.6642	&	0.2358	&	 0.6396	&	 0.2238	 &	 0.6168	&	0.2129	&	0.5953	&	0.2031	&	 0.5753	 \\					
10	&&	0.3097	&	0.7356	&	0.2917	&	0.7083	&	0.2758	&	 0.6829	&	 0.2615	 &	 0.6589	&	0.2486	&	0.6365	&	0.2370	&	 0.6154	 \\					
11	&&	    	&	    	&	    	&	    	&	0.3171	&	 0.7242	&	 0.3004	 &	 0.6996	&	0.2854	&	0.6761	&	0.2719	&	 0.6542	 \\					
12	&&	    	&	    	&	    	&	    	&	    	&	    	 &	    	 &	    	 &	0.3239	&	0.7146	&	0.3083	&	 0.6917	 \\ &&&&&&&&&&&&&\\ \hline					
$x$	&&	\multicolumn{2}{c}{$n=25$}			&	 \multicolumn{2}{c}{$n=26$}			 &	 \multicolumn{2}{c}{$n=27$}			 &	 \multicolumn{2}{c}{$n=28$}			&	 \multicolumn{2}{c}{$n=29$}			 &	\multicolumn{2}{c}{$n=30$}			\\ \hline					 
0	&&	0.0000	&	0.1078	&	0.0000	&	0.1039	&	0.0000	&	 0.1002	&	 0.0000	 &	 0.0969	&	0.0000	&	0.0937	&	0.0000	&	 0.0907	 \\					
1	&&	0.0020	&	0.1775	&	0.0020	&	0.1713	&	0.0019	&	 0.1654	&	 0.0018	 &	 0.1599	&	0.0018	&	0.1548	&	0.0017	&	 0.1500	 \\					
2	&&	0.0162	&	0.2363	&	0.0156	&	0.2280	&	0.0150	&	 0.2204	&	 0.0145	 &	 0.2132	&	0.0140	&	0.2064	&	0.0135	&	 0.2001	 \\					
3	&&	0.0333	&	0.2895	&	0.0320	&	0.2795	&	0.0308	&	 0.2703	&	 0.0296	 &	 0.2616	&	0.0286	&	0.2534	&	0.0276	&	 0.2457	 \\					
4	&&	0.0556	&	0.3392	&	0.0534	&	0.3277	&	0.0513	&	 0.3169	&	 0.0494	 &	 0.3068	&	0.0476	&	0.2974	&	0.0460	&	 0.2885	 \\					
5	&&	0.0793	&	0.3862	&	0.0761	&	0.3733	&	0.0731	&	 0.3612	&	 0.0704	 &	 0.3498	&	0.0678	&	0.3392	&	0.0654	&	 0.3291	 \\					
6	&&	0.1059	&	0.4312	&	0.1016	&	0.4169	&	0.0976	&	 0.4036	&	 0.0939	 &	 0.3910	&	0.0905	&	0.3792	&	0.0873	&	 0.3681	 \\					
7	&&	0.1336	&	0.4744	&	0.1281	&	0.4589	&	0.1230	&	 0.4444	&	 0.1183	 &	 0.4308	&	0.1140	&	0.4179	&	0.1099	&	 0.4058	 \\					
8	&&	0.1634	&	0.5161	&	0.1565	&	0.4995	&	0.1503	&	 0.4839	&	 0.1445	 &	 0.4692	&	0.1391	&	0.4554	&	0.1341	&	 0.4423	 \\					
9	&&	0.1941	&	0.5565	&	0.1859	&	0.5388	&	0.1784	&	 0.5222	&	 0.1714	 &	 0.5065	&	0.1650	&	0.4917	&	0.1590	&	 0.4777	 \\					
10	&&	0.2264	&	0.5956	&	0.2167	&	0.5769	&	0.2078	&	 0.5594	&	 0.1997	 &	 0.5428	&	0.1921	&	0.5271	&	0.1851	&	 0.5123	 \\					
11	&&	0.2596	&	0.6335	&	0.2484	&	0.6140	&	0.2381	&	 0.5955	&	 0.2287	 &	 0.5781	&	0.2200	&	0.5616	&	0.2119	&	 0.5460	 \\					
12	&&	0.2942	&	0.6703	&	0.2813	&	0.6499	&	0.2696	&	 0.6307	&	 0.2588	 &	 0.6124	&	0.2488	&	0.5952	&	0.2396	&	 0.5788	 \\					
13	&&	0.3297	&	0.7058	&	0.3152	&	0.6848	&	0.3018	&	 0.6648	&	 0.2896	 &	 0.6460	&	0.2784	&	0.6280	&	0.2680	&	 0.6110	 \\					
14	&&	    	&	    	&	    	&	    	&	0.3352	&	 0.6982	&	 0.3214	 &	 0.6786	&	0.3088	&	0.6600	&	0.2972	&	 0.6423	 \\					
15	&&	    	&	    	&	    	&	    	&	    	&	    	 &	    	 &	    	 &	0.3400	&	0.6912	&	0.3271	&	 0.6729	 \\ \hline					
\end{tabular}

\newpage

\normalsize
\noindent
Table A2:  Limits of 99\% equal-tail OLC interval when $x$ successes are observed in a sample of size $n$


\scriptsize

\noindent
\begin{tabular}{cccccccccccccc} \hline


$x$	&&	\multicolumn{2}{c}{ $n=1$}		&	\multicolumn{2}{c}{$n=2$}			 &	 \multicolumn{2}{c}{$n=3$}			&	 \multicolumn{2}{c}{$n=4$}			 &	 \multicolumn{2}{c}{$n=5$}			&	\multicolumn{2}{c}{$n=6$}			 \\  \hline					 
0	&&	0.0000	&	0.9900	&	0.0000	&	0.8801	&	0.0000	&	 0.7572	&	 0.0000	 &	 0.6565	&	0.0000	&	0.5761	&	0.0000	&	 0.5121	 \\					
1	&&	0.0100	&	1.0000	&	0.0050	&	0.9950	&	0.0033	&	 0.9297	&	 0.0025	 &	 0.8441	&	0.0020	&	0.7636	&	0.0017	&	 0.6928	 \\					
2	&&		&		&		&		&	0.0703	&	0.9967	&	0.0500	 &	0.9500	 &	 0.0388	&	0.8841	&	0.0318	&	0.8177	\\					 
3	&&		&		&		&		&		&		&		&		&	 0.1159	&	 0.9612	 &	 0.0925	&	0.9075	\\ 					

&&&&&&&&&&&&&\\ \hline					
$x$	&&	\multicolumn{2}{c}{$n=7$}			&	 \multicolumn{2}{c}{$n=8$}			 &	 \multicolumn{2}{c}{$n=9$}			 &	 \multicolumn{2}{c}{$n=10$}			&	 \multicolumn{2}{c}{$n=11$}			 &	\multicolumn{2}{c}{$n=12$}			\\ \hline					 
0	&&	0.0000	&	0.4602	&	0.0000	&	0.4175	&	0.0000	&	 0.3819	&	 0.0000	 &	 0.3518	&	0.0000	&	0.3260	&	0.0000	&	 0.3037	 \\					
1	&&	0.0014	&	0.6320	&	0.0013	&	0.5800	&	0.0011	&	 0.5353	&	 0.0010	 &	 0.4966	&	0.0009	&	0.4629	&	0.0008	&	 0.4334	 \\					
2	&&	0.0269	&	0.7558	&	0.0233	&	0.7005	&	0.0206	&	 0.6513	&	 0.0184	 &	 0.6078	&	0.0167	&	0.5693	&	0.0152	&	 0.5351	 \\					
3	&&	0.0770	&	0.8510	&	0.0660	&	0.7965	&	0.0577	&	 0.7462	&	 0.0513	 &	 0.7004	&	0.0462	&	0.6590	&	0.0420	&	 0.6217	 \\					
4	&&	0.1490	&	0.9230	&	0.1261	&	0.8739	&	0.1094	&	 0.8253	&	 0.0966	 &	 0.7794	&	0.0866	&	0.7368	&	0.0784	&	 0.6976	 \\					
5	&&		&		&		&		&	0.1747	&	0.8906	&	0.1532	 &	0.8468	 &	 0.1365	&	0.8046	&	0.1232	&	0.7649	\\					 
6	&&		&		&		&		&		&		&		&		&	 0.1954	&	 0.8635	 &	 0.1754	&	0.8246	\\ 				

&&&&&&&&&&&&&\\ \hline					
$x$	&&	\multicolumn{2}{c}{$n=13$}			&	 \multicolumn{2}{c}{$n=14$}			 &	 \multicolumn{2}{c}{$n=15$}			 &	 \multicolumn{2}{c}{$n=16$}			&	 \multicolumn{2}{c}{$n=17$}			 &	\multicolumn{2}{c}{$n=18$}			\\ \hline					 
0	&&	0.0000	&	0.2842	&	0.0000	&	0.2670	&	0.0000	&	 0.2518	&	 0.0000	 &	 0.2382	&	0.0000	&	0.2260	&	0.0000	&	 0.2150	 \\					
1	&&	0.0008	&	0.4072	&	0.0007	&	0.3840	&	0.0007	&	 0.3633	&	 0.0006	 &	 0.3446	&	0.0006	&	0.3277	&	0.0006	&	 0.3124	 \\					
2	&&	0.0140	&	0.5045	&	0.0130	&	0.4771	&	0.0121	&	 0.4524	&	 0.0113	 &	 0.4300	&	0.0106	&	0.4097	&	0.0100	&	 0.3912	 \\					
3	&&	0.0385	&	0.5879	&	0.0356	&	0.5574	&	0.0331	&	 0.5297	&	 0.0309	 &	 0.5045	&	0.0289	&	0.4815	&	0.0272	&	 0.4604	 \\					
4	&&	0.0717	&	0.6618	&	0.0660	&	0.6291	&	0.0612	&	 0.5991	&	 0.0570	 &	 0.5717	&	0.0534	&	0.5464	&	0.0502	&	 0.5233	 \\					
5	&&	0.1122	&	0.7280	&	0.1030	&	0.6938	&	0.0953	&	 0.6622	&	 0.0886	 &	 0.6331	&	0.0828	&	0.6061	&	0.0777	&	 0.5812	 \\					
6	&&	0.1592	&	0.7875	&	0.1458	&	0.7527	&	0.1345	&	 0.7200	&	 0.1248	 &	 0.6897	&	0.1165	&	0.6614	&	0.1092	&	 0.6351	 \\					
7	&&	0.2125	&	0.8408	&	0.1939	&	0.8061	&	0.1785	&	 0.7731	&	 0.1653	 &	 0.7420	&	0.1539	&	0.7128	&	0.1441	&	 0.6855	 \\					
8	&&		&		&		&		&	0.2269	&	0.8215	&	0.2097	 &	0.7903	 &	 0.1949	&	0.7606	&	0.1822	&	0.7326	\\					 
9	&&		&		&		&		&		&		&		&		&	 0.2394	&	 0.8051	 &	 0.2233	&	0.7767	\\ 					

&&&&&&&&&&&&&\\ \hline					
$x$	&&	\multicolumn{2}{c}{$n=19$}			&	 \multicolumn{2}{c}{$n=20$}			 &	 \multicolumn{2}{c}{$n=21$}			 &	 \multicolumn{2}{c}{$n=22$}			&	 \multicolumn{2}{c}{$n=23$}			 &	\multicolumn{2}{c}{$n=24$}			\\ \hline					 
0	&&	0.0000	&	0.2049	&	0.0000	&	0.1958	&	0.0000	&	 0.1875	&	 0.0000	 &	 0.1798	&	0.0000	&	0.1728	&	0.0000	&	 0.1662	 \\					
1	&&	0.0005	&	0.2984	&	0.0005	&	0.2857	&	0.0005	&	 0.2739	&	 0.0005	 &	 0.2631	&	0.0004	&	0.2531	&	0.0004	&	 0.2438	 \\					
2	&&	0.0095	&	0.3743	&	0.0090	&	0.3587	&	0.0086	&	 0.3444	&	 0.0082	 &	 0.3311	&	0.0078	&	0.3188	&	0.0075	&	 0.3074	 \\					
3	&&	0.0257	&	0.4411	&	0.0244	&	0.4232	&	0.0232	&	 0.4067	&	 0.0221	 &	 0.3915	&	0.0211	&	0.3773	&	0.0202	&	 0.3640	 \\					
4	&&	0.0473	&	0.5019	&	0.0448	&	0.4821	&	0.0425	&	 0.4637	&	 0.0405	 &	 0.4467	&	0.0386	&	0.4308	&	0.0369	&	 0.4160	 \\					
5	&&	0.0733	&	0.5581	&	0.0693	&	0.5367	&	0.0657	&	 0.5167	&	 0.0624	 &	 0.4982	&	0.0595	&	0.4809	&	0.0569	&	 0.4647	 \\					
6	&&	0.1028	&	0.6106	&	0.0970	&	0.5878	&	0.0919	&	 0.5665	&	 0.0873	 &	 0.5466	&	0.0832	&	0.5280	&	0.0794	&	 0.5106	 \\					
7	&&	0.1354	&	0.6599	&	0.1278	&	0.6359	&	0.1209	&	 0.6135	&	 0.1148	 &	 0.5925	&	0.1092	&	0.5727	&	0.1042	&	 0.5542	 \\					
8	&&	0.1710	&	0.7062	&	0.1611	&	0.6814	&	0.1523	&	 0.6580	&	 0.1444	 &	 0.6360	&	0.1374	&	0.6153	&	0.1309	&	 0.5958	 \\					
9	&&	0.2093	&	0.7498	&	0.1969	&	0.7243	&	0.1860	&	 0.7002	&	 0.1762	 &	 0.6774	&	0.1674	&	0.6559	&	0.1595	&	 0.6356	 \\					
10	&&	0.2502	&	0.7907	&	0.2351	&	0.7649	&	0.2218	&	 0.7403	&	 0.2100	 &	 0.7169	&	0.1993	&	0.6947	&	0.1897	&	 0.6737	 \\					
11	&&		&		&		&		&	0.2597	&	0.7782	&	0.2456	 &	0.7544	 &	 0.2329	&	0.7317	&	0.2215	&	0.7102	\\					 
12	&&		&		&		&		&		&		&		&		&	 0.2683	&	 0.7671	 &	 0.2549	&	0.7451	\\ 					

&&&&&&&&&&&&&\\ \hline					
$x$	&&	\multicolumn{2}{c}{$n=25$}			&	 \multicolumn{2}{c}{$n=26$}			 &	 \multicolumn{2}{c}{$n=27$}			 &	 \multicolumn{2}{c}{$n=28$}			&	 \multicolumn{2}{c}{$n=29$}			 &	\multicolumn{2}{c}{$n=30$}			\\ \hline					 
0	&&	0.0000	&	0.1602	&	0.0000	&	0.1545	&	0.0000	&	 0.1493	&	 0.0000	 &	 0.1444	&	0.0000	&	0.1398	&	0.0000	&	 0.1355	 \\					
1	&&	0.0004	&	0.2352	&	0.0004	&	0.2271	&	0.0004	&	 0.2196	&	 0.0004	 &	 0.2126	&	0.0003	&	0.2060	&	0.0003	&	 0.1998	 \\					
2	&&	0.0072	&	0.2968	&	0.0069	&	0.2869	&	0.0066	&	 0.2776	&	 0.0064	 &	 0.2689	&	0.0062	&	0.2607	&	0.0059	&	 0.2530	 \\					
3	&&	0.0193	&	0.3517	&	0.0185	&	0.3401	&	0.0178	&	 0.3293	&	 0.0172	 &	 0.3192	&	0.0166	&	0.3096	&	0.0160	&	 0.3006	 \\					
4	&&	0.0353	&	0.4022	&	0.0339	&	0.3892	&	0.0326	&	 0.3771	&	 0.0314	 &	 0.3656	&	0.0302	&	0.3548	&	0.0292	&	 0.3446	 \\					
5	&&	0.0544	&	0.4495	&	0.0522	&	0.4352	&	0.0501	&	 0.4218	&	 0.0482	 &	 0.4092	&	0.0465	&	0.3973	&	0.0448	&	 0.3861	 \\					
6	&&	0.0760	&	0.4942	&	0.0728	&	0.4788	&	0.0699	&	 0.4643	&	 0.0672	 &	 0.4506	&	0.0647	&	0.4377	&	0.0624	&	 0.4255	 \\					
7	&&	0.0996	&	0.5367	&	0.0954	&	0.5203	&	0.0916	&	 0.5048	&	 0.0880	 &	 0.4902	&	0.0847	&	0.4763	&	0.0817	&	 0.4632	 \\					
8	&&	0.1251	&	0.5774	&	0.1198	&	0.5600	&	0.1149	&	 0.5436	&	 0.1104	 &	 0.5281	&	0.1062	&	0.5134	&	0.1023	&	 0.4995	 \\					
9	&&	0.1523	&	0.6163	&	0.1457	&	0.5982	&	0.1397	&	 0.5810	&	 0.1341	 &	 0.5646	&	0.1290	&	0.5492	&	0.1242	&	 0.5345	 \\					
10	&&	0.1810	&	0.6537	&	0.1731	&	0.6348	&	0.1658	&	 0.6169	&	 0.1592	 &	 0.5999	&	0.1530	&	0.5837	&	0.1473	&	 0.5683	 \\					
11	&&	0.2112	&	0.6896	&	0.2018	&	0.6701	&	0.1933	&	 0.6515	&	 0.1854	 &	 0.6339	&	0.1782	&	0.6171	&	0.1715	&	 0.6010	 \\					
12	&&	0.2429	&	0.7241	&	0.2319	&	0.7040	&	0.2219	&	 0.6849	&	 0.2128	 &	 0.6667	&	0.2044	&	0.6493	&	0.1966	&	 0.6328	 \\					
13	&&	0.2759	&	0.7571	&	0.2633	&	0.7367	&	0.2518	&	 0.7172	&	 0.2413	 &	 0.6985	&	0.2316	&	0.6806	&	0.2227	&	 0.6635	 \\					
14	&&		&		&		&		&	0.2828	&	0.7482	&	0.2709	 &	0.7291	 &	 0.2599	&	0.7108	&	0.2498	&	0.6933	\\					 
15	&&		&		&		&		&		&		&		&		&	 0.2892	&	 0.7401	 &	 0.2778	&	0.7222	\\ 				

\hline					
\end{tabular}

\end{appendix}



\newpage
\vspace{0.1in}
\section*{References}
\begin{description}
\item Agresti, A. and Coull, B. A. (1998). Approximate is better than ``exact'' for interval estimation of binomial proportions. {\em The American Statistician,} 52, 119--126.
\item Agresti, A. and Gottard, A. (2005). Comment: Randomized confidence intervals and the mid-$P$ approach. {\em Statistical Science,} 20, 367--371.
\item Agresti, A. and Gottard, A. (2007). Nonconservative exact small-sample inference for discrete data. {\em Computational Statistics and Data Analysis,} 51, 6447--6458.
\item Berry, G. and Armitage, P. (1995). Mid-$p$ confidence intervals: a brief review. {\em The Statistician,} 44, 417--423.
\item Blyth, C. R. and Still, H. A. (1983). Binomial confidence intervals. {\em Journal of the American Statistical Association,} 78, 209--212.
\item Brown, L. D., Cai, T. T. and DasGupta, A. (2001). Interval estimation for a binomial proportion (with discussion). {\em Statistical Science,} 16, 108--116.
\item Casella, G. (1986). Refining binomial confidence intervals. {\em Canadian Journal of Statistics,} 14, 113--129.
\item Crow, E. L.  (1956). Confidence intervals for a proportion. {\em Biometrika,} 43, 423--435.
\item Ghalanos, A. and Theussl, S. (2015). Package Rsolnp. CRAN Repository.
\item Jovanovic, B. D. and Levy, P. S. (1997). A look at the rule of three. {\em The American Statistician,} 51, 137--139.
\item Leemis, L. M. and Trivedi, K. S. (1996). A comparison of approximate interval estimators for the Bernoulli Parameter. {\em The American Statistician,} 50, 63--68.
\item Mehta, c. R. and Walsh, S. J. (1992). Comparison of exact, Mid-$p$, and Mantel-H\'{a}enszel confidence intervals for the common odds ratio across several $2 \times 2$ contingency tables. {\em The American Statistician,} 46, 146--150.
\item Newcombe, R. G. (1998). Two-sided confidence intervals for the single proportion: comparison of seven methods. {\em Statistics in Medicine,} 17, 857--872.
\item Schilling, M. F. and Doi, J. A. (2014). A coverage probability approach to finding an optimal binomial confidence procedure. {\em The American Statistician,} 68, 133--145.
\item Vollset, S. E. (1993). Confidence intervals for a binomial proportion. {\em Statistics in Medicine,} 12, 809--824.
\item Wilson, E. B. (1927). Probable inference, the law of succession, and statistical inference. {\em Journal of the American Statistical Association,} 22, 209--212.
\end{description}

\end{document}


\begin{frontmatter}
\title{Locally correct confidence intervals for a binomial proportion: A new criteria for an
interval estimator\\ (Supplementary material)}
\runtitle{Locally correct confidence intervals (Supplementary material)}

\begin{aug}
\author{\fnms{Paul H.} \snm{Garthwaite}},
\author{\fnms{Maha W.} \snm{Moustafa}}
\and
\author{\fnms{Fadlalla G.} \snm{Elfadaly}}
\address{School of Mathematics and Statistics, The Open University, Milton Keynes, MK7 6AA, UK.}
\runauthor{Garthwaite et al.}
\end{aug}

\begin{keyword}[class=MSC2020]
\kwd[Primary ]{62A01}
\kwd{62F25}
\end{keyword}

\begin{keyword}
\kwd{Clopper-Pearson}
\kwd{coverage}
\kwd{discrete distribution}
\kwd{mid-$p$}
\kwd{shortest interval}
\end{keyword}

\end{frontmatter}


\noindent
Proofs of Propositions 1--5 are given below.\vspace{.12in}

\setcounter{equation}{0}
\renewcommand{\theequation}{A.\arabic{equation}}

\noindent
{\bf Proof of Proposition 1}\vspace{.05in}

\noindent
Let $\beta =1-\alpha$. The values of $\beta$ that are of interest in Proposition 1 are those in the range $(0.73, \,  0.9999)$. The plan is to take tiny slices that together cover this range. One value of $n$ is considered at a time. When the endpoints of a slice are $r$ and $s$ $(r<s)$, let $v_0, \ldots , v_{n-1}$ be the values of $u_0, \ldots , u_{n-1}$ when $\beta=r$ and let $w_0, \ldots , w_{n-1}$ be their values when $\beta=s$. For the slices that we consider, we will show that $P(X \geq i |p=v_i) \geq s$ for $i=0,\ldots ,n-1$ and that each of $\delta u_{n-1} / \delta \beta, \ldots, \delta u_{0} / \delta \beta$ is positive for $r \leq \beta \leq s$. This is sufficient to prove Proposition 1, as it gives the monotonicity property and also establishes that $P(X \geq i |p= \beta) \geq \beta$ for $r \leq \beta \leq s$.

Lemma A1 gives a set of slices for which each slice satisfies  $P(X \geq i |p=v_i) \geq s$ for $i=0,\ldots ,n-1$. Lemma A2 gives a formula for $\delta u_{n-1} / \delta \beta$ and shows that it is a positive non-decreasing function of $\beta$ for $0<\beta<1$ and any $n$. Lemma A3 gives minimum and maximum limits for $\delta u_{i-1} / \delta \beta$ within a slice in terms of minimum and maximum limits for $\delta u_{i} / \delta \beta$ in that slice. Using an iterative procedure, direct computation for a slice can establish that $\delta u_{n-1} / \delta \beta, \ldots, \delta u_{0} / \delta \beta$ are positive for values of $\beta$ in that slice. The result in Lemma A2 is used to start an iteration and Lemma A3 gives the iterative step.

Repetitive computation yielded the result in the first lemma.\vspace{.1in}\\
{\bf Lemma A1.} \, Partition the interval $(0.73, \,  0.9999)$ into non-overlapping touching sub-intervals (slices) of width 0.0001 in the range 0.73--0.95, width 0.00005 in the range 0.95--0.99, width 0.00001 in the range 0.99--0.999, and width 0.000001 in the range 0.999--0.9999. Let $(r,s)$ denote any of these slices. Then for any $n$ in $1,\ldots ,200$,
\begin{equation} \label{eqNA1}
P(X \geq i |p=v_i) \geq s
\end{equation}
for $i=0,\ldots ,n-1$, where  $v_0, \ldots , v_{n-1}$ are the values of $u_0, \ldots , u_{n-1}$ when $\beta=r$.\vspace{.1in}\\
{\bf Lemma A2.} \, If
 $u_{n-1}$ is the upper limit of a upper one-sided interval with confidence level $\beta$ when $X=n-1$, then
\begin{equation} \label{eqNA2}
    \frac{d u_{n-1}}{d \beta} = \frac{(n+1)(1-u_{n-1})^2}{1 -(n+1)u_{n-1}^n + n u_{n-1}^{n+1}}.
\end{equation}
Also $d u_{n-1} /d \beta$ is a monotonic decreasing function of $\beta$ and $d u_{n-1} /d \beta >0$ for $0<\beta < 1$.\vspace{.1in}\\
{\bf Proof.} \begin{equation} \label{eqNA3}
    \beta = \frac{1}{1-u_{n-1}} \int_{u_{n-1}}^1 p^n \, dp
\end{equation}
so
\begin{equation} \label{eqNA4}
   \frac{d\beta}{du_{n-1}} = \frac{1}{(1-u_{n-1})^2} \int_{u_{n-1}}^1 p^n \, dp \, - \, \frac{u_{n-1}^n}{1-u_{n-1}} = \frac{\beta - u_{n-1}^n}{1-u_{n-1}}.
\end{equation}
From equation (\ref{eqNA3}), $\beta > u_{n-1}^n$, so equation (\ref{eqNA4}) implies that $d\beta /du_{n-1} >0$ (and also $ du_{n-1} /d\beta >0$).\\
Equation (\ref{eqNA3}) also gives
\[
    \beta  =\frac{1-u_{n-1}^{n+1}}{(n+1)(1-u_{n-1})}.
\]
Hence,
\[
    \frac{d\beta}{du_{n-1}} = \frac{1 -(n+1)u_{n-1}^n + n u_{n-1}^{n+1}}{(n+1)(1-u_{n-1})^2}
\]
and equation (\ref{eqNA2}) follows.
Further differentiation yields
\[
   (n+1) \frac{d^2\beta}{du_{n-1}^2} = \frac{2 \{ 1 -(n+1)u_{n-1}^n + n u_{n-1}^{n+1} \} }{(1-u_{n-1})^3}
   - \frac{n(n+1)(u_{n-1}^{n-1} -u_{n-1}^n) }{(1-u_{n-1})^2},
\]
giving
\[
   (n+1) \frac{d^2\beta}{du_{n-1}^2} = \frac{ A }{(1-u_{n-1})^3}
\]
where
\[
A =  2 - n(n+1)u_{n-1}^{n-1} + 2(n^2 -1) u_{n-1}^{n} - (n^2-n) u_{n-1}^{n+1}.
\]
Clearly $d^2\beta / du_{n-1}^2 >0$ if and only if $A>0$. Now
\begin{eqnarray*}
    \frac{dA}{du_{n-1}} & = & -n(n+1)(n-1)u_{n-1}^{n-2} + 2n(n^2 -1) u_{n-1}^{n-1} - (n+1)(n^2-n) u_{n-1}^{n} \\
    & = & -n(n^2-1)(1-u_{n-1})^2u_{n-1}^{n-1}.
\end{eqnarray*}
Hence $dA / du_{n-1}$ is negative for $u_{n-1} <1$. Also, $A=0$ when $u_{n-1} =1$. Hence $A$ is positive for $u_{n-1} <1$. Consequently, for $u_{n-1} <1$, we have that $d^2\beta / du_{n-1}^2 >0$ and $d\beta / du_{n-1}$ is a monotonic increasing function of $u_{n-1}$. As $d\beta / du_{n-1} >0$, it follows that $ du_{n-1} / d\beta$ is a monotonic decreasing function of $\beta$. \hfill $\Box$
\vspace{.15in}

\noindent
{\bf Lemma A3.} \,  Let $v_0, \ldots , v_{n-1}$ and $w_0, \ldots , w_{n-1}$ denote the values of $u_0, \ldots , u_{n-1}$ when $\beta=r$ and $\beta=s$, respectively, where $r<s$. Let $(\delta u_i / \delta \beta )_{\mbox{min}}$ and $(\delta u_i / \delta \beta )_{\mbox{max}}$ be lower and upper bounds on $(\delta u_i / \delta \beta )$ for $\beta \in [r,s]$ and suppose $(\delta u_i / \delta \beta )_{\mbox{min}} >0$. Then \vspace{.08in}
\begin{equation} \label{eqNA5}
\frac{v_i-w_{i-1}- \{ P(X \geq i \, | \, p=w_{i}) - r \} \mbox{$(\delta u_i / \delta \beta )_{\mbox{max}}$} }{ s- P(X \geq i \, | \, p=v_{i-1})}
\end{equation}
and
\begin{equation} \label{eqNA6}
\frac{w_i-v_{i-1}- \{ P(X \geq i \, | \, p=v_{i}) - s \} \mbox{$(\delta u_i / \delta \beta )_{\mbox{min}}$} }{ r- P(X \geq i \, | \, p=w_{i-1})}
\end{equation}
are lower and upper bounds on $(\delta u_{i-1} / \delta \beta )$ for $\beta \in [r,s]$ if the lower bound in (\ref{eqNA5}) is positive.\vspace{.1in}\\
{\bf Proof.} \, Differentiation of
\[
    \beta = \frac{1}{u_i  -u_{i-1}} \int_{u_{i-1}}^{u_i} \, P(X \geq i \, | \, p) \, dp.
\]
gives
\begin{eqnarray} \label{eqNA7}
    \delta \beta & = & - \frac{\delta u_i - \delta u_{i-1}}{(u_i - u_{i-1})^2} \int_{u_{i-1}}^{u_i} \, P(X \geq i \, | \, p) \, dp +
    \frac{\delta u_i}{u_i - u_{i-1}}  P(X \geq i \, | \, p=u_i) \nonumber \\
    & & - \, \frac{\delta u_{i-1}}{u_i - u_{i-1}}  P(X \geq i \, | \, p=u_{i-1}).
\end{eqnarray}
Consequently,
\begin{eqnarray}
    (u_i -u_{i-1} ) \delta \beta & = & - (\delta u_i - \delta u_{i-1}) \beta  +
    \delta u_i  P(X \geq i \, | \, p=u_i) \nonumber \\
    & & - \, \delta u_{i-1}  P(X \geq i \, | \, p=u_{i-1}) \nonumber
\end{eqnarray}
\[
    = \delta u_i \{ P(X \geq i \, | \, p=u_i) - \beta \}
    + \delta u_{i-1} \{ \beta  -  P(X \geq i \, | \, p=u_{i-1}) \}.
\]
Hence,
\begin{equation} \label{eqNA8}
    \frac{ \delta u_{i-1}}{\delta \beta} \, = \, \left[ (u_i -u_{i-1} ) - \{ P(X \geq i \, | \, p=u_i) - \beta \}  \frac{ \delta u_{i}}{\delta \beta} \right] \{ \beta  -  P(X \geq i \, | \, p=u_{i-1}) \}^{-1}.
\end{equation}
To place bounds on $\delta u_{i-1} / \delta \beta $ we form bounds on the terms in the right-hand side of equation(\ref{eqNA8}). If $u_{i-1} \in [v_{i-1}, \, w_{i-1}]$, then lower and upper bounds on $u_i -u_{i-1}$ are $v_i -w_{i-1}$ and $w_i -v_{i-1}$, respectively. Lower and upper bounds on $\{ P(X \geq i \, | \, p=u_{i}) - \beta \}$ are $\{ P(X \geq i \, | \, p=v_{i}) - s \} $ and $\{ P(X \geq i \, | \, p=w_{i} ) - r \} $, respectively. Also, lower and upper bounds on $\{ \beta-P(X \geq i \, | \, p=u_{i-1}) \}$ are $\{r- P(X \geq i \, | \, p=w_{i-1}) \} $ and $\{s- P(X \geq i \, | \, p=v_{i-1} \} $, respectively. Lastly, by assumption $(\delta u_i / \delta \beta )_{\mbox{min}}$ and $(\delta u_i / \delta \beta )_{\mbox{max}}$ are lower and upper bounds on $(\delta u_i / \delta \beta )$ for $\beta \in [r,s]$. All these bounds are positive, either automatically or by assumption. As long as $u_{i-1} \in [v_{i-1}, \, w_{i-1}]$,  they can replace the terms in the right-hand side of equation(\ref{eqNA8}) to give (\ref{eqNA5}) and (\ref{eqNA6}) as bounds on $ \delta u_{i-1} / \delta \beta$ for $\beta \in [r,s]$.

To prove Lemma A3, it remains to show that $u_{i-1} \in [v_{i-1}, \, w_{i-1}]$ if $\beta \in [r,s]$. Let $\epsilon = v_i-w_{i-1}- \{ P(X \geq i \, | \, p=w_{i}) - r \} (\delta u_i / \delta \beta )_{\mbox{max}}$ and suppose the quantity in (\ref{eqNA5}) is positive. Then $\epsilon > 0$ and, from equation(\ref{eqNA8}), $\delta u_i / \delta \beta  > \epsilon /2$ if $u_{i-1} <w_{i-1} + \epsilon /2$ and $\beta \in [r,s]$. Thus $u_{i-1}$ would be an increasing function of $\beta$ if it was in the interval $(w_{i-1}, \, w_{i-1} + \epsilon/2)$. However, clearly $u_{i-1}$ is a continuous function of $\beta$ and $u_{i-1} = w_{i-1}$ when $\beta=s$. It follows that $u_{i-1}$ cannot be greater than $w_{i-1}$ for $\beta <s$ when the assumptions in Lemma A3 hold. Also, $\delta u_{i-1} / \delta \beta$ is positive for $u_{i-1} \in (v_{i-1}, \, w_{i-1})$ if $\beta \in [r,s]$ and $u_{i-1}=v_{i-1}$ when $\beta=r$. Hence $u_{i-1} \geq v_{i-1}$ for $\beta \in [r,s]$.\hfill $\Box$
\vspace{.15in}

The interval $(0.73, \, 0.9999)$ was partitioned into slices exactly as in Lemma 1. Each slice and each $n$ $(n=1, \ldots,200)$ was taken in turn. Bounds on $\delta u_{n-1} / \delta \beta$ were calculated using Lemma A2 and it was checked that the lower bound was positive. Then bounds on $\delta u_{n-2} / \delta \beta, \ldots, \delta u_{o} / \delta \beta$ were calculated sequentially using Lemma A3. With the set of slices given in Lemma A1, the lower bounds were all positive so the assumptions made in Lemma A3 were always satisfied, which completes the proof of Proposition 1.\vspace{.15in}

\noindent
{\bf Proof of Proposition 2}\vspace{.05in}

\noindent
For $x=0, \ldots , n$, the $1-\alpha$ upper-tail interval estimate for $p$ given by the OLC method is $(0,u_x^*)$ when $X=x$ and $X \sim\mbox{bin}(n, \, p)$. Suppose $p \leq u_x^*$. We have that $0 < u_0^* < \ldots < u_{x-1}^* < u_{x}^* < \ldots < u_n^* =1$, so $p$ is contained in the confidence interval if $X \geq x$. If $p$ is equally likely to take any value in the interval $(0, \, 1)$, then $\mbox{P}(X = x) =1/(n+1)$ for $x=0, \ldots ,n$. Hence, the AEL of the OLC method is $ \sum_{x=0}^{n} u_x^* /(n+1)$. We aim to give conditions on $n$ and $\alpha$ under which this AEL is less than that of any other method that yields $1-\alpha$ upper-tail LCC intervals.

We will refer to $(u_0, \ldots, u_n)$ as the {\em partition} of an interval estimator when $(0,u_x)$ is the upper-tail interval given by the estimator for a $1-\alpha$ confidence level when $X=x$. Suppose that the method yielding $1-\alpha$ upper-tail LCC intervals with the smallest AEL gives the partition $(b_0, \ldots, b_n)$. Let $g_n(x, p) = \mbox{P}(X \geq x \, | \, p)$. Since the method is locally correct,
\begin{equation} \label{eqA4}
\frac{1}{b_{x}-b_{x-1}} \int_{b_{x-1}}^{b_x } g_n(x, p) \, dp \geq 1-\alpha
\end{equation}
for $x=1, \ldots ,n$, and $0 < b_0 \leq \ldots \leq b_x \leq b_{x+1} \leq \ldots \leq b_n =1.$

In Lemma A7 we will determine conditions on the $u_x^*$ under which
\begin{equation} \label{eqA6}
b_x=u_x^* \quad  \mbox{for } x=0, \ldots, n.
\end{equation}
The conditions that are obtained can be examined for individual values of $n$ and $\alpha$. We then prove Proposition 2 by showing that the conditions hold for $1 \leq n \leq 200$ and $0.0001 < \alpha <0.27 $.

Lemma A7 is proved by induction, by showing that if $b_x = u_x^*$ for $x=n, n-1, \ldots , k+1$ ($k \geq 0$), then it also holds for $x=k$. To start the induction, note that $1=b_n = u_n^* $. Lemmas A4 and A5 give preparatory results that place bounds on the $b_x$.\vspace{.1in}\\
{\bf Lemma A4}. \,  Suppose $b_{x+1} = u_{x+1}^*$. Then (i) $b_{x} \geq u_{x}^*$ if $x \leq n-1$; (ii) $b_{x-1}  \leq u_{x-1}^*$ if $x \geq 1$;  (iii)  $b_0 =u_0^*$ if $x=1$; and (iv) $b_{x-1}+b_x \leq u_{x-1}^*+u_x^* $ if $x \geq 1$.\vspace{.1in}\\
{\bf Proof.} \,  As $g_n(x+1, p)$ is a monotonic strictly decreasing function of $p$, $\int_{c}^{u_{x+1}^*} g_n(x+1, p) \, dp  / (u_{x}^*-c )  <  \int_{u_{x}^*}^{u_{x+1}^*} g_n(x+1, p) \, dp  / (u_{x+1}^*-u_{x}^*)$ for any $c < u_x^*$. Given $u_{x+1}^*$, the OLC method chooses $u_x^*$ as the value for which $\int_{u_{x}^*}^{u_{x+1}^*} g_n(x+1, p) \, dp  / (u_{x+1}^*-u_{x}^*) = 1- \alpha$. Thus $\int_{c}^{u_{x+1}^*} g_n(x+1, p) \, dp  / (u_{x+1}^*-c) < 1-\alpha$ for any $c<u_x^*$. By assumption, $b_{x+1}=u_{x+1}^*$, so (\ref{eqA4}) gives result (i), that $b_x \geq u_x^*$. For (ii), note that $(b_0, \ldots, b_{x-1}, u_x^*, \ldots , u_n)$ would be a partition that gave LCC intervals if $b_{x-1} \geq u_{x-1}^*$. Since $\sum_{i=0}^{x-1} b_i  + \sum_{i=x}^{n} u_i^* = \sum_{i=0}^{n} b_i +(u_x^* -b_x)$ and (from (i)) $b_x \geq u_x^*$, this would contradict the assumption that $(b_0, \ldots, b_n)$ is the partition of minimum average length that satisfies (\ref{eqA4}), unless $b_{x-1} = u_{x-1}^*$. Consequently, $b_{x-1} \leq u_{x-1}^*$ if $b_{x+1} = u_{x+1}^*$. For (iii), suppose that $b_1 =u_1^*$. Then the partition $(u_0^*, b_1, b_2, \ldots, b_n)$ gives LCC intervals. As $(b_0, \ldots, b_n)$ is the partition giving LCC intervals that has the shortest average length, it follows that $b_0 \leq u_0^*$ so, from (i), $b_0 = u_0^*$. For (iv), from (ii) we have that $b_{x-1} \leq u_{x-1}^*$, so the average coverage for $p \in (b_{x-2}, \, u_{x-1}^*)$ is not less than the average coverage for $p \in (b_{x-2}, \, b_{x-1})$. Hence, as $(b_0 , \ldots , b_n)$ is a partition that gives LCC intervals, so does the partition $(b_0 , \ldots , b_{x-2}, u_{x-1}^*, u_x, b_{x+1}, \ldots,  b_n)$. Now $(b_0 , \ldots , b_n)$ is the partition that gives the {\em shortest} average length of LCC interval estimators, so we have that $b_{x-1}+b_{x} \leq u_{x-1}^*+u_{x}^*$. \quad \quad { \large $\diamond$}
\vspace{.2in}

To place a further bound on $b_{x-1}$, let $p_{x-1}^{\#}$ be the value of $p$ for which $g_n(x-1,p)= 1-\alpha$. Analogous to equation (\ref{eqA4}) we have
\begin{equation} \label{eqA2b}
\frac{1}{b_{x-1}-b_{x-2}} \int_{b_{x-2}}^{b_{x-1 }} g_n(x-1, p) \, dp \geq 1-\alpha,
\end{equation}
so $b_{x-1} \geq p_{x-1}^{\#}$, since $g_n(x-1,p)$ is a monotonically increasing function of $p$. In combination with Lemma A1, this gives the following result.\vspace{.1in}\\
{\bf Lemma A5.} \, If $b_{x+1}=u_{x+1}^*$ and $x \geq 1$, then
\begin{equation} \label{eqAn1}
u_x^* \leq b_x \leq \min(u_{x+1}^*, u_{x-1}^*+u_x^* - p_{x-1}^{\#})
\end{equation}
where $p_{x-1}^{\#}$ is given by $g_n(x-1, p_{x-1}^{\#})= 1-\alpha$.
\vspace{.2in}\\
{\bf Proof.} \, For the lower limit, $u_x^* \leq b_x$ from part (i) of Lemma A4. For the upper limit, first $b_x \leq b_{x+1}$ so $b_x \leq u_{x+1}^*$. Also, from part (iv) of Lemma A4, $b_x \leq u_{x-1}^*+u_x^* -b_{x-1}$, so $b_x \leq u_{x-1}^*+u_x^* -p_{x-1}^{\#}$. \quad \quad { \large $\diamond$}

\vspace{.2in}

We will require  a quantity $A_n(x, p^*,y)$ that is defined by:
\begin{equation} \label{eqAn2}
A_n(x, p^*,y) = \frac{1}{2y} \int_{p^*-y}^{p^*+y} g_n(x, p) \, dp \quad \quad \mbox{for } \, 0 \leq y \leq \min(p^*, \,1-p^*).
\end{equation}
Thus  $A_n(x, p^*,y)$ is the average coverage given by $g_n(i, p) $ in the interval $(p^*-y, \, p^*+y)$. Its differential with respect to $y$ is denoted as $A_n'(x, p^*,y)$. The following lemma gives an important characteristic of $A_n(x, p^*,y)$ that underlies conditions for the OLC method to yield LCC intervals with minimum average length.\vspace{.1in}\\
{\bf Lemma A6.} \, Suppose that $0 < y_1 <y_2  \leq \min(p^*, \,1-p^*)$ and that $A_n(x, p^*,y_1) >A_n(x, p^*,y_2)$. Then,
\begin{equation} \label{eqAn3}
A_n(x, p^*,y_1) >A_n(x, p^*,y) \quad \quad \mbox{for any } y \in (y_1, \, y_2 ]
\end{equation}
if $A_n'(x, p^*,y_1) <0$.
\vspace{.2in}\\
{\bf Proof.} \,  Differentiation of $A_n(x, p^*,y)$ with respect to $y$ gives
\begin{equation} \label{eqAn4}
A_n'(x, p^*,y) = - \frac{1}{2y^2} \int_{p^*-y}^{p^*+y} g_n(x, p) \, dp +\frac{1}{2y}  \{ g_n(x, p^*+y) +g_n(x, p^*-y)  \} .
\end{equation}
For sufficiently small $y$, the smooth function $g_n(x,p)$ is approximately linear for $p$ in the interval $(p^*-y,\, p^*+y)$, so $A_n(x, p^*,y)$ does not vary as $y \rightarrow 0$, giving $A_n'(x, p^*,0)=0$.

Let $Q_n(x, p^*,y) = 2y^2 \,A_n'(x, p^*,y)$. Then the differential of $Q_n(x, p^*,y)$ with respect to $y$ is:
\begin{eqnarray}
Q_n'(x, p^*,y) &=& -  \{ g_n(x, p^*+y) +g_n(x, p^*-y)  \}  + \{ g_n(x, p^*+y) +g_n(x, p^*-y)  \} \nonumber \\
& &+ \, y \, \{ g_n'(x, p^*+y) -g_n'(x, p^*-y)  \}  \nonumber \\
&=& y \, \{ g_n'(x, p^*+y) -g_n'(x, p^*-y) \}. \label{eqnA9}
\end{eqnarray}
Put $h(y) = g_n'(x, p^*+y) /g_n'(x, p^*-y)$. The differential of $h(y)$ with respect to $y$ is:
\begin{equation} \label{eqAn5}
h'(y)= \psi(y) [ c_1 - c_2 y^2] ,
\end{equation}
where
\[ \psi(y) = 2 (p^*+y)^{x-2} (1-p^* -y)^{n-x-1} (p^*-y)^{-x} (1-p^* +y)^{-n+x-1}. \]
\[ c_1 = p^* (1-p^*) \{ x-1-p^* (n-1) \}, \]
and
\[ c_2 = p^* (n-1)-n+x   .\]
The behavior of $A_n'(x, p^*,y)$ depends upon the signs of $c_1$ and $c_2$.\vspace{.2in}

Case 1: $c_1 \leq 0$; $c_2>0$ \, or \, $c_1<0$ ; $c_2=0$.\\
As $\psi >0$ for all $y$, the following results follow sequentially.
\begin{enumerate}
	\item[(i)] $h'(y) <0$ for any $y$, so $h(y)$ is a monotonic decreasing function of $y$.
	\item[(ii)] As $h(0)=1$, it follows that $h(y)<1$ for all $y>0$, so $g_n'(x, p^*+y) -g_n'(x, p^*-y)$ is negative for any $y>0$.
	\item[(iii)] $Q_n'(x, p^*,y)$ is negative for   $y>0$, so $Q_n(x, p^*,y)$ is a monotonic decreasing function of $y$.
	\item[(iv)] As $Q_n(x, p^*,0)=0$,  it follows that $Q_n(x, p^*,0)$ is negative for any  $y>0$. Thus $A_n'(x, p^*,y)$ is also negative for any  $y>0$.
	\item[(v)] $A_n(x, p^*,y)$ is a monotonic decreasing function of $y$, so $A_n(x, p^*,y_1)>A_n(x, p^*,y)$ for any $y> y_1$. Thus Lemma A6 holds if $c_1 \leq 0$ and $c_2>0$, or if $c_1<0$ and $c_2=0$.
\end{enumerate}
\vspace{.2in}

Case 2: $c_1 \geq 0$; $c_2 \leq0$.\\
Now $h'(y) \geq 0$ for any $y$ and reasoning similar to that for Case 1 shows that $A_n(x, p^*,y_1) \leq A_n(x, p^*,y)$ for any $y> y_1$. Thus the conditions required by Lemma A6 do not hold, so there is nothing to verify.\vspace{.2in}

Case 3: $c_1 <0$; $c_2<0$.\\
The following results hold.
\begin{enumerate}
	\item[(i)] $h'(y) <0$ for $y<(c_1/c_2)^{1/2}$ and $h'(y) >0$ for $y>(c_1/c_2)^{1/2}$, so $h(y)$ is a $\cup \,$-shaped function of $y$.
	\item[(ii)] Suppose $h(y)$ is below 1 until $y=y^*$. As $h(0)=1$, it follows that $g_n'(x, p^*+y) -g_n'(x, p^*-y)$ is negative for $0<y<y^*$ and positive for   $y>y^*$.
	\item[(iii)] $Q_n'(x, p^*,y)$ is negative for $0<y<y^*$ and positive for   $y>y^*$, so $Q_n(x, p^*,y)$ is a $\cup \,$-shaped function of $y$.
	\item[(iv)] Suppose $Q_n(x, p^*,y)$ is below 0 until $y=y^{\#}$.
	As $Q_n(x, p^*,0) =0$ and $Q_n(x, p^*,y) = 2y^2A_n'(x, p^*,y)$,  it follows that both $Q_n(x, p^*,y)$ and $A_n'(x, p^*,y)$  are negative for $0<y<y^{\#}$ and positive for  $y>y^{\#}$. Hence, $A_n(x, p^*,y)$ is a $\cup \,$-shaped function of $y$.
	\item[(v)] If $y_2 >y_1$ and $A_n(x, p^*,y_1) >A_n(x, p^*,y_2)$, then $A_n(x, p^*,y_1)>A_n(x, p^*,y)$ for any $y \in (y_1, \, y_2)$. Thus Lemma A6 holds if $c_1 < 0$ and $c_2<0$.
	\vspace{.2in}
\end{enumerate}
Notes. \, In (i) it is assumed that $(c_1/c_2)^{1/2}  \leq \min(p^*, \,1-p^*)$. If $(c_1/c_2)^{1/2} > \min(p^*, \,1-p^*)$, then $h'(y) <0$ for all feasible values of $y$ and Case 1 applies. Similarly, in (ii), if $h(y)$ never reaches 1, then $Q_n'(x, p^*,y)$ is always negative and Case 1 again applies (c.f. Case 1, parts (iii)-(v)). Likewise, in (iv), if  $Q_n(x, p^*,y)$ is always below 0 for any feasible values of $y$, then Case 1 applies (parts (iv) and (v)). \vspace{.2in}

Case 4: $c_1 > 0$; $c_2>0$.\\
The following results hold.
\begin{enumerate}
	\item[(i)] $h'(y) >0$ for $y<(c_1/c_2)^{1/2}$ and $h'(y) <0$ for $y>(c_1/c_2)^{1/2}$, so $h(y)$ is a $\cap \,$-shaped function of $y$.
	\item[(ii)] Suppose $h(y)$ drops below 1 at $y=y^*$. As $h(0)=1$, it follows that $g_n'(x, p^*+y) -g_n'(x, p^*-y)$ is positive for $0<y<y^*$ and negative for   $y>y^*$.
	\item[(iii)] $Q_n'(x, p^*,y)$ is positive for $0<y<y^*$ and negative for   $y>y^*$, so $Q_n(x, p^*,y)$ is a $\cap \,$-shaped function of $y$.
	\item[(iv)] Suppose $Q_n(x, p^*,y)$ drops below 0 at $y=y^{\#}$. As $Q_n(x, p^*,0) =0$,  it follows that both $Q_n(x, p^*,y)$ and $A_n'(x, p^*,y)$ are positive for $0<y<y^{\#}$ and negative for $y>y^{\#}$. Thus, $A_n(x, p^*,y)$ is a $\cap \,$-shaped function of $y$.
	\item[(v)] If $A_n'(x, p^*,y_1)<0$, then $A_n(x, p^*,y_1)>A_n(x, p^*,y)$ for any $y> y_1$. Thus Lemma A6 holds if $c_1 > 0$ and $c_2>0$.
	\vspace{.2in}
\end{enumerate}
Notes. \, In (i) it is assumed that $(c_1/c_2)^{1/2}  \leq \min(p^*, \,1-p^*)$. If $(c_1/c_2)^{1/2} > \min(p^*, \,1-p^*)$, then $h'(y) >0$ for all feasible values of $y$ and Case 2 applies. Similarly, Case 2 applies if $h(y)$ never drops below 1 in (ii) for any feasible value of $y$, or if  $Q_n(x, p^*,y)$ never drops below 0 in (iii). \vspace{.2in}\\
This completes the proof of Lemma A6, as Cases 1 - 4 cover all combinations of $c_1$ and $c_2$. \quad \quad { \large $\diamond$}
\vspace{.2in}\\
{\bf Lemma A7.} \,For $x=1, \ldots, n-1$, let $p_x^* = (u_x^* +u_{x-1}^*)/2$ and let $\eta_x = u_x^* -p_x^*$. Also, define $p_{x-1}^{\#}$ by $g_n(x-1,p_{x-1}^{\#} )= 1-\alpha$ for $x=2, \ldots, n-1$ and $p_{0}^{\#}=0$. Let $\xi_x =\min(u_{x+1}^*- p_x^*, \, p_x^*-p_{x-1}^{\#})$.  Suppose
\begin{equation} \label{eqA10}
A_n(x,p_x^*,\xi_x) < 1-\alpha
\end{equation}
and
\begin{equation} \label{eqA11}
A_n'(x,p_x^*,\eta_x)<0
\end{equation}
both hold for $x=1, \ldots, n-1$.\,
Then $b_x =u_x^*$ for $x=0, \ldots,n$ and $(u_0^*, \ldots, u_n^*)$ is the partition that yields locally correct $1-\alpha$ confidence intervals of shortest average length.\vspace{.1in}\vspace{.2in}\\
{\bf Proof.} \,  We have that $b_n=u_n^*$ and, from Lemma A4, $b_0 =u_0^*$ if $b_1 = u_1^*$. Hence we must show that $b_x =u_x^*$ for $x=1,\ldots , n-1$.
We will prove the result by induction: we assume that  $b_x =u_x^*$ for
$x=n, n-1, \ldots , k+1$ ($k \geq 1$) and will show this implies that $b_k =u_k^*$.

Now $A_n(k,p_k^*,\eta_k) =  (2 \eta_k)^{-1} \int_{p_k^* -\eta_k} ^{p_k^*+\eta_k} g_n(k,p) \, dp = (2 \eta_k)^{-1} \int_{u_{k-1}} ^{u_k} g_n(k,p) \, dp = 1-\alpha$. Hence, by assumption, $A_n(k,p_k^*,\eta_k) > A_n(k,p_k^*,\xi_k)$. Also $A_n'(k,p_k^*,\eta_k)<0$ so, from Lemma A6, $A_n(k,p_k^*,\eta_k) > A_n(k,p_k^*,y)$ for any $y \in (\eta_k, \, \xi_k ] $. Consequently,
\begin{equation} \label{eqA16}
1-\alpha > A_n(k,p_k^*,y) \quad \quad \mbox{for any } y \in (\eta_k, \, \xi_k].
\end{equation}
Put $b_k =p_k^* +\tau_k$. From part (i) of Lemma A4,  $\tau_k \geq 0$, and from part (iv), $b_{k-1} \leq u_{k-1}^*+u_k^* -b_k \, = \, 2p_k^* -b_k = p_k^* -\tau _k$. Hence, as $g_n(k,p)$ is a monotonic strictly increasing function of $p$, we have
\[
1-\alpha \leq \frac{1}{b_k-b_{k-1}} \int_{b_{k-1}} ^{b_k} g_n(k,p) \, dp \leq \frac{1}{2\tau_k} \int_{p_k^*-\tau_k} ^{p_k^*+\tau_k} g_n(k,p) \, dp = A_n(k,p_k^*,\tau_k).
\]


From equation (\ref{eqA16}), it follows that $\tau_k$ is not within the interval $(\eta_k, \, \xi_k ] $. However, $u_k^*= p_k^* +\eta_k$ so, from Lemma A5,
\[ p_k^* +\eta_k \leq b_k \leq \min(u_{k+1}^*,  \,u_{k-1}^*+u_k^*-p_{k-1}^{\#}) =p_k^*+\xi_k. \]
Thus $\eta_k \leq \tau_k \leq \xi_k$. It follows that $\tau_k =\eta_k$, so $b_k = u_k^*$. \quad \quad {\large $\diamond$}\vspace{.2in}

The following lemma extends Lemma A7 to a range of values of $1-\alpha$. 
\vspace{.15in}

\noindent
{\bf Lemma A8.} \,  Let $v_0, \ldots , v_{n}$ and $w_0, \ldots , w_{n}$ denote the values of $u_0^*, \ldots , u_{n}^*$ when $1-\alpha=r$ and $1-\alpha=s$, respectively, where $r<s$. Define $p_x^*$, $\eta_x$ and $\xi_x$ as in Lemma A7 and define $p_{x-1}^{\min}$ and $p_{x-1}^{\max}$ by $g_n(x-1,p_{x-1}^{\min})=r$ and $g_n(x-1,p_{x-1}^{\max})=s$ for $x=2,\ldots,n$. Let $p_{0}^{\min} = p_{0}^{\max}=0$ and define $\xi_x^{\min}$ and $\xi_x^{\max}$ by $\xi_x^{\min} = \min \{ v_{x+1}-(w_x +w_{x-1})/2, \, (v_x+v_{x-1})/2 - p_{x-1}^{\max} \} $ and $\xi_x^{\max} = \min \{ w_{x+1}-(v_x +v_{x-1})/2, \, (w_x+w_{x-1})/2 - p_{x-1}^{\min} \} $ for $x=1,\ldots, n-1$. Then, for $r \leq 1-\alpha < s$ and $x=1,\ldots, n-1$,
\begin{enumerate}
\item[(a)]  $A_n'(x,p_x^*,\eta_x)$ is negative  if
\begin{equation} \label{NA20}
2r > g_n(x,w_x) + g_n(x,w_{x-1}).
\end{equation}
\item[(b)] $A_n(x,p_x^*,\xi_x) < 1-\alpha$ if
\begin{equation} \label{NA21}
r > (\xi_x^{\max} +\xi_x^{\min})^{-1} \int_{(w_x+w_{x-1})/2-\xi_x^{\min}}^{(w_x+w_{x-1})/2+\xi_x^{\max}} g_n(x, p) \, dp.
\end{equation}
\end{enumerate}
{\bf Proof of (a).} \, From equation (\ref{eqAn4}),
\[ A_n'(x,p_x^*,\eta_x) =
 - \frac{1}{2\eta_x^2} \int_{p_x^*-\eta_x}^{p_x^*+\eta_x} g_n(x, p) \, dp +\frac{1}{2\eta_x}  \{ g_n(x, p_x^*+\eta_x) +g_n(x, p_x^*-\eta_x)  \} .
\]
Now $p_x^*+\eta_x=u_x^*$, $p_x^*-\eta_x=u_{x-1}^*$ and $2\eta_x = u_x^* - u_{x-1}^*$, so
\[ \frac{1}{2\eta_x} \int_{p_x^*-\eta_x}^{p_x^*+\eta_x} g_n(x, p) \, dp =
\frac{1}{ u_x^* - u_{x-1}^*} \int_{u_{x-1}^*}^{u_x^*} g_n(x, p) \, dp =1-\alpha. \]
Also, $g_n(x, p_x^*+\eta_x)= g_n(x, u_x^*) $ and $g_n(x, p_x^*-\eta_x)= g_n(x, u_{x-1}^*)$. Consequently, $\eta_x A_n'(x,p_x^*,\eta_x) =-(1-\alpha) +  \{ g_n(x, u_x^*) +g_n(x, u_{x-1}^*)  \} /2 $. When $r \leq 1-\alpha < s$, the maximum values of $g_n(x, u_x^*)$ and $g_n(x, u_{x-1}^*)$ are $g_n(x, w_x)$ and $g_n(x, w_{x-1})$, while the minimum value of $1-\alpha$ is $r$. Hence $A_n'(x,p_x^*,\eta_x)$ is negative for $r \leq 1-\alpha <s$ if $2r > g_n(x,w_x) + g_n(x,w_{x-1})$.\vspace{.01in}\\
{\bf Proof of (b).} \, From equation (\ref{eqAn2})
\begin{equation} \label{NA22}
A_n(x, p_x^*,\xi_x) = \frac{1}{2\xi_x} \int_{p_x^*-\xi_x}^{p_x^*+\xi_x} g_n(x, p) \, dp
\end{equation}
Thus $A_n(x, p_x^*,\xi_x)$ is the average value of $g_n(x, p)$ between two limits. Since $g_n(x, p)$ is a monotonic increasing function of $p$, the average value of $g_n(x, p)$ between limits $a$ and $b$ will increase as $a$ increases and as $b$ increases. For $r \leq 1-\alpha <s$, an upper bound for $p_x^*+\xi_x$ is $(w_x+w_{x-1})/2+\xi_x^{\max}$ and an upper bound for $p_x^*-\xi_x$ is $(w_x+w_{x-1})/2-\xi_x^{\min}$. Hence an upper bound on $A_n(x, p_x^*,\xi_x)$ is $(\xi_x^{\max} +\xi_x^{\min})^{-1} \int_{(w_x+w_{x-1})/2-\xi_x^{\min}}^{(w_x+w_{x-1})/2+\xi_x^{\max}} g_n(x, p) \, dp$ \, for $r \leq 1-\alpha <s$. Also, for these values of $1-\alpha$ a lower bound on $1-\alpha$ is $r$.
It follows that $A_n(x,p_x^*,\xi_x) < 1-\alpha$ if
$r > (\xi_x^{\max} +\xi_x^{\min})^{-1} \int_{(w_x+w_{x-1})/2-\xi_x^{\min}}^{(w_x+w_{x-1})/2+\xi_x^{\max}} g_n(x, p) \, dp$. \quad \quad {\large $\diamond$}\vspace{.2in}

Proposition 2 relates to values of $1-\alpha$ in the range $(0.73,0.9999)$. As in the proof of Proposition 1, tiny slices were taken that together cover this range and one value of $n$ was taken at a time. Specifically, the interval $(0.73,0.9999)$ was partitioned into non-overlapping touching slices of width 0.001 in range 0.73-0.99, width 0.0001 in the range 0.99-0.999, and width 0.00001 in the range 0.999-0.9999. Each $n$ was taken in turn $(n=2, \ldots, 200)$ and for each slice it was checked that the inequalities in (\ref{NA20}) and (\ref{NA21}) held, where $r$ and $s$ are the endpoints of a slice. This verifies that $A_n'(x,p_x^*,\eta_x) <0$ and  $A_n(x,p_x^*,\xi_x) < 1-\alpha$ for $0.73<1-\alpha <0.9999$. Consequently, the conditions in (\ref{eqA10}) and (\ref{eqA11}) hold for these values of $1-\alpha$, so Proposition 2 follows from Lemma A7.\vspace{.15in}

\noindent
{\bf Proof of Proposition 3}\vspace{.08in}\\
As in the proof of Proposition 1, the interval $(0.73, \,  0.9999)$ was partitioned into non-overlapping touching sub-intervals (slices) of width 0.0001 in the range 0.73--0.95, width 0.00005 in the range 0.95--0.99, width 0.00001 in the range 0.99--0.999, and width 0.000001 in the range 0.999--0.9999. Each slice was taken in turn and it was checked that $L(x,n+1,s) <L(x,n,r)$ and $U(x,n+1,s) <U(x,n,r)$ for $n=1,\ldots, 200$, $x=0, \ldots, n$, where $r$ and $s$ are the endpoints of the slice ($r<s$). From the monotonicity property in Proposition 1, it follows that $L(x,n+1,\alpha) <L(x,n,\alpha)$ and $U(x,n+1,\alpha) <U(x,n,\alpha)$ for $r \leq \alpha <s$, and together the slices cover the range $0.0001 \leq \alpha < 0.27$, as required.\quad \quad {\large $\diamond$}\vspace{.2in}

\noindent
{\bf Proof of Proposition 4}\vspace{.08in}\\
For each value of $n$ in turn ($1 \leq n \leq 200 $), the following algorithm was run for the upper limits.
\begin{enumerate}
\item Set $\alpha^* =0.1$.
\item Using the mid-$p$ method, determine the upper endpoints (spikes) for $1-\alpha^*$ upper-tail confidence intervals for $X=0,1, \ldots, n$.
\item Determine the minimum average coverage between consecutive spikes and denote this by $\beta $. In the initial loop (when $\alpha^* =0.1$) check that $\beta$ is greater than $1-\alpha^*$.
\item Stop if $\beta> 0.9999 $. Otherwise, set $\alpha^*$ equal to $1-\beta$ and return to step 2.
\end{enumerate}
For each  value of $n$ the check was satisfied that $\beta > 1-\alpha^*$ for $\alpha^* =0.1$. Also, for each $n$ the value of $\beta$ increased in each cycle of the loop until $\beta$ exceeded 0.9999. (The increase in $\beta$ was sometimes very small so average coverages were calculated to an accuracy of 1.0E$-10$.)
Now, as $\alpha $ decreases, the endpoints given by the mid-$p$ method for a $1- \alpha$ upper-tail interval increase, so the minimum average coverage between spikes also increases. Hence for $1-\alpha^* < 1-\alpha \leq \beta$, the minimum coverage between spikes of $1-\alpha$ upper-tail intervals exceeds the minimum coverage between spikes of $1-\alpha^*$ upper-tail intervals. This latter coverage is $\beta$. Hence for $1-\alpha^* < 1-\alpha \leq \beta$ the mid-$p$ method gives LCC intervals. For each $n$ the set of loops start with $\alpha^* =0.1$ and end with $\beta > 0.9999$, so the mid-$p$ method gives LCC intervals for $0.0001 \leq \alpha \leq 0.1$. Symmetry implies that the equivalent result holds for lower-tail intervals. \quad \quad {\large $\diamond$}\vspace{.2in}

\noindent
{\bf Proof of Proposition 5}\vspace{.08in}\\
Suppose $p=p^*$ and $0 \leq u_0 \leq \ldots \leq u_n =1$.
If $1- \textstyle{\frac{1}{2}} \mbox{P}(X = 0 \, | \, p=p^*) \leq 1-\alpha $, put $i^*=0$. Otherwise, define $i^*$ by
\begin{equation} \label{B17}
\mbox{P}(X  \geq i^* \, | \, p=p^*)  - \textstyle{\frac{1}{2}} \mbox{P}(X = i^* \, | \, p=p^*)
\leq 1-\alpha
\end{equation}
and
\begin{equation} \label{B18}
\mbox{P}(X \geq i^*-1 \, | \, p=p^*)  - \textstyle{\frac{1}{2}} \mbox{P}(X = i^*-1 \, | \, p=p^*) > 1-\alpha.
\end{equation}
The coverage of a method of forming confidence intervals must equal $\mbox{P}(X \geq i \, | \, p=p^*)$ for some $i$.
The value midway between $\mbox{P}(X \geq i-1)$ and $\mbox{P}(X \geq i)$ is $\mbox{P}(X \geq i-1) - \frac{1}{2}\mbox{P}(X = i-1)$. Similarly, $\mbox{P}(X \geq i) -\frac{1}{2}\mbox{P}(X = i)$ is midway between $\mbox{P}(X \geq i)$ and $\mbox{P}(X \geq i+1)$. Hence, if $i^* \geq 1$, the feasible coverage that is closest to equalling $1-\alpha$ is $\mbox{P}(x \geq i^* \, | \, p=p^*)$. A method of forming confidence intervals achieves this coverage if $u_{i^*-1} < p^* \leq u_{i^*}$. If $i^*=0$, the coverage closest to $1-\alpha$ is 1, which is the coverage when $p^* \leq u_0$.


Let $\tilde{u}_i$ denote the upper limit given by the mid-$p$ method when $X=i$. From equation~(13), $\tilde{u}_i $ satisfies
\begin{equation} \label{B19}
\mbox{P}(X \geq i \, | \, p=\tilde{u}_i)  - \textstyle{\frac{1}{2}} \mbox{P}(X = i \, | \, p=\tilde{u}_i)
 \,  =  \, 1-\alpha ,
\end{equation}
for $i \leq n-1$. Suppose $i^* \leq n-1$.
Putting $i=i^*$ in (\ref{B19}) and comparison with (\ref{B17}) yields
\begin{equation} \label{B20}
\mbox{P}(X \geq i^*  \, | \, p=p^*)  - \textstyle{\frac{1}{2}} \mbox{P}(X = i^* \, | \, p=p^*)
\, \leq \,
\mbox{P}(X \geq i^*  \, | \, p=\tilde{u}_{i^*})  - \textstyle{\frac{1}{2}} \mbox{P}(X = i^* \, | \, p=\tilde{u}_{i^*}),
\end{equation}
so $ p^* \leq \tilde{u}_{i^*}$.
Hence the mid-$p$ method has the feasible coverage closest to $1-\alpha$ when $i^* =0$. As $\tilde{u}_n =1$, we also have $ p^* \leq \tilde{u}_{i^*}$ for $i^*=n$.
When $i^* \geq 1$, putting $i=i^*-1$ in (\ref{B19}) and comparison with (\ref{B18}) similarly yields
\[
\mbox{P}(X \geq i^*-1 \, | \, p=p^*)  - \textstyle{\frac{1}{2}} \mbox{P}(X = i^*-1 \, | \, p=p^*) \makebox[1in]{ }
\]
\begin{equation} \label{B21}
\makebox[1in]{ }
 \,  > \,
\mbox{P}(X \geq i^*-1  \, | \, p=\tilde{u}_{i^*-1})  - \textstyle{\frac{1}{2}} \mbox{P}(X = i^*-1 \, | \, p=\tilde{u}_{i^*-1}),
\end{equation}
so $ p^* > \tilde{u}_{i^*-1}$. Thus, when $1 \leq i^* \leq n$, we have that $ \tilde{u}_{i^*-1}< p^* \leq \tilde{u}_{i^*}$. Hence, for any $p^*$, the coverage of the mid-$p$ method is as close to $1-\alpha$ as the coverage of any method that meets the conditions of Proposition 5.\quad \quad { \large $\diamond$}\vspace{.25in}